\documentclass[a4paper,11pt,oneside]{article}
\usepackage{amsmath,amsthm,amsfonts,amssymb}
\usepackage{a4}
\usepackage[ansinew]{inputenc}
\usepackage{mathrsfs}

\newcommand{\shO}{\mathcal{O}}
\newcommand{\shC}{\mathscr{C}}
\newcommand{\shE}{\mathcal{E}}
\newcommand{\shT}{\mathcal{T}}

\newcommand{\R}{\mathbf{R}}
\newcommand{\C}{\mathbf{C}}
\newcommand{\Z}{\mathbf{Z}}

\newcommand{\End}{\rm{End }}
\newcommand{\Der}{\rm{Der}}

\newcommand{\Ad}{{\rm{Ad }}}
\newcommand{\Sp}{{\rm{Sp }}}
\newcommand{\ad}{{\rm{ad}}}

\newcommand{\id}{{\rm id}}
\newcommand{\tr}{{\rm{tr}}}
\newcommand{\str}{{\rm{str}}}
\newcommand{\SL}{{\rm SL}}
\newcommand{\rmS}{{\rm S}}
\newcommand{\rmO}{{\rm O}}
\newcommand{\rmU}{{\rm U}}
\newcommand{\SO}{{\rm SO}}
\newcommand{\GL}{{\rm GL}}
\newcommand{\OSp}{{\rm OSp}}
\newcommand{\SOSp}{{\rm SOSp}}
\newcommand{\PSL}{{\rm PSL}}
\newcommand{\Aut}{{\rm Aut}}

\newcommand{\osp}{{\mathfrak{osp}}}
\newcommand{\gl}{{\mathfrak{gl}}}
\newcommand{\psl}{{\mathfrak{psl}}}

\newcommand{\red}{\rm{red}}

\newtheorem{Thm}{Theorem}[section]
\newtheorem{Def}[Thm]{Definition}
\newtheorem{Prop}[Thm]{Proposition}
\newtheorem{Cor}[Thm]{Corollary}
\newtheorem{Lemma}[Thm]{Lemma}

\theoremstyle{remark}
\newtheorem*{Rem}{Remark}

\numberwithin{equation}{section}

\begin{document}

%\thispagestyle{empty}
%\begin{center}\vspace*{2 cm}{\huge \bf Riemannian Supergeometry} \\[3 cm]
%Inaugural-Dissertation\\
%zur\\
%Erlangung des Doktorgrades\\
%der Mathematisch-Naturwissenschaftlichen Fakult\"at\\
%der Universit\"at zu K\"oln\\[5 cm]
%vorgelegt von\\
%{\bf Oliver Goertsches} \\
%aus K\"oln\\[3 cm]
%K\"oln\\[1 cm]
%2006
%\end{center}

%\newpage
%\thispagestyle{empty}
%\vspace*{1 pt}
%\vfill \noindent
%\begin{tabular}{ll}
%Berichterstatter: & Prof. Dr. Gudlaugur Thorbergsson\\
%& Prof. Dr. Uwe Semmelmann\\
%\\
%Tag der m\"undlichen Pr\"ufung: & im Sommer 2006
%\end{tabular} 

%\pagestyle{empty}
\title{Riemannian Supergeometry}
\author{O. Goertsches\footnote{Funded by DFG, SFB TR/12}}
%\begin{center}
%\huge Riemannian Supergeometry\\[0.7 cm]
%\large O. Goertsches\footnote{Funded by DFG, SFB TR/12}\\[1 cm]
%\end{center}
\maketitle

%\thispagestyle{empty}
%\noindent
%\centerline{\bf Abstract}\\

\begin{abstract} Motivated by Zirnbauer \cite{Zirnbauer}, we develop a theory of Riemannian supermanifolds up to a definition of Riemannian symmetric superspaces. Va\-rious fundamental concepts needed for the study of these spaces both from the Riemannian and the Lie theoretical viewpoint are introduced, e.g.~geodesics, isometry groups and invariant metrics on Lie supergroups and homogeneous superspaces.\end{abstract}

\section{Introduction}

Although there exists a theory of differential geometry of supermanifolds\footnote{There exist various versions of supergeometry; the formalism we use is the sheaf-theoretic approach by Berezin, Kostant and Leites \cite{Berezin}, \cite{Kostant}, \cite{Leites}.} -- see Deligne and Morgan \cite{DelMor}, Schmitt \cite{Schmitt} or Varadarajan \cite{Varadara} -- a notion of Riemannian metric for supermanifolds only seldom occurs in the literature. More precisely, beyond the existence of a Levi-Civita connection -- see Monterde and S\'anchez-Valenzuela \cite{MontSan} -- no general theory of Riemannian supermanifolds is available.

The motivation for the development of such a theory came from the physicists: in 1996, Zirnbauer \cite{Zirnbauer} defined Riemannian symmetric superspaces to be a quotient of complex Lie supergroups, together with some distinguished Riemannian symmetric space embedded into the underlying manifold. 
%%%%%%%%%%%% IM ARTIKEL WEG
%Note that, with this definition, the question whether a homogeneous superspace is Riemannian symmetric is decided only in the underlying manifold. 

%In the standard theory, the fact that Riemannian symmetric spaces may be regarded as homogeneous spaces is very relevant; nevertheless, the Riemannian viewpoint is essential for a mathematical treatment of these spaces to such an extent that most textbooks define them as Riemannian manifolds with the symmetry property of the existence of geodesic symmetries. In spirit of this, our aim is to give a similar definition of Riemannian symmetric superspaces and afterwards recognize them as special homogeneous superspaces. But before that, various fundamental concepts have to be introduced and studied, e.g.~geodesics, isometry groups and invariant metrics on Lie supergroups and homogeneous superspaces.
%%%%%%%%%%%%% IM ARTIKEL DAZU

Our aim is to give a definition of these objects similar to the standard theory, namely as Riemannian supermanifolds with a symmetry property, and afterwards recognize them as special homogeneous superspaces. But before that, various fundamental concepts have to be introduced and studied, e.g.~geodesics, isometry groups and invariant metrics on Lie supergroups and homogeneous superspaces.

It is to be mentioned that the non-linear theory developed in this thesis already has some infinitesimal counterpart in the mathematical literature; for example, Cort\'es \cite{Cortes} defines a notion of infinitesimal pseudo-Riemannian symmetric superspace, and Serganova \cite{Serganova} lists up involutive automorphisms of the simple Lie superalgebras over $\R$ and $\C$. 

In the sections 2 and 3 we review the basics of differential supergeometry and Lie supergroups; we assume some familiarity with linear superalgebra, see e.g.~Varadarajan \cite{Varadara} or Deligne and Morgan \cite{DelMor}.

{\emph{Acknowledgements.}} This paper essentially is my thesis which I wrote under the supervision of G. Thorbergsson at the University of Cologne. He deserves my gratitude for the constantly encouraging support during the last years.

\section{Foundations of Supergeometry}

\subsection{Supermanifolds and their Morphisms}\label{Section_Smfd}

The model in the category of supermanifolds is the space $\R^{n\mid m}$, which is by definition the ringed space consisting of the topological space $\R^n$ and the sheaf of super $\R$-algebras %\footnote{As mentioned in the introduction, we assume familiarity with the basic notions of linear superalgebra, see e.g.~Varadarajan \cite{Varadara} or Deligne and Morgan \cite{DelMor}.} 
$\shC^\infty\otimes\Lambda_\R[\xi_1,\ldots,\xi_m]$.

A {\it supermanifold (graded manifold)} of dimension $n|m$ is a ringed space \mbox{$M=(|M|,\shO_M)$}, where $|M|$ is a topological space (Hausdorff, countable base) and the structural sheaf $\shO_M$ is a sheaf of super $\R$-algebras with unity, locally isomorphic to $\R^{n\mid m}$. Sections of the structural sheaf are referred to as {\it superfunctions} on $M$; if there is no danger of confusion, we simply call them {\it functions}. 

Let $M$ and $N$ be supermanifolds. A {\it morphism} $\Phi:M\to N$ is a morphism of ringed spaces: $\Phi=(\phi,\phi^*)$, where $\phi$ is a continuous map between the underlying topological spaces and $\phi^*:\shO_N\to \phi_*\shO_M$ is a morphism of sheaves of \mbox{$\R$-super} algebras with unity. The morphism is {\it not} determined by the map between the topological spaces; nevertheless, the map on global sections $\phi^*_N:\shO_N(N)\to \phi_*\shO_M(N)=\shO_M(M)$ determines the whole morphism, i.e.~$\phi$ and $\phi^*$, cf. \cite{Kostant}, p. 208.

The nilpotent functions define an ideal sheaf $J$ of $\shO_M$, and the ringed space $M_{\red}:=(|M|,\shO_M/J)$ is a differentiable manifold. We call it the {\it reduced} or {\it underlying manifold} or the {\it support} of $M$. The quotient map $\shO_M\to \shO_M/J$ defines a morphism $M_{\red}\to M$ sending a superfunction $f$ on $M$ to a smooth function $\widetilde{f}$ on the reduced manifold. The {\it value} of a superfunction $f$ at some point $p\in M_{\red}$ is defined to be $\widetilde{f}(p)$, which coincides with the unique real number $\lambda$ such that $f-\lambda$ is not invertible as an element of the stalk $\shO_{M,p}$; sometimes, we simply write $f(p)$ for the value of $f$ of $p$ although this might be misleading -- since $f$ is not determined by all its values, associating to $f$ the function sending $p$ to $f(p)$ is not an injective mapping. Hence, this does not provide us with a realization of $\shO_M$ as a sheaf of ordinary real-valued functions. 

If $U\subset M_{\red}$ is such that $\shO_M(U)= \shC^\infty(U)\otimes \Lambda_\R[\xi_1,\ldots,\xi_m]$ and coordinates $x_i$ of the reduced manifold on $U$ are given, then we call $(x_i,\xi_\alpha)$ {\it coordinates of} $M$ on $U$. Note that our convention is to give roman indices for even (here: the {\it even coordinates} $x_i$) and greek indices for odd objects (here: the {\it odd coordinates} $\xi_\alpha$) -- this will be convenient for notation purposes. However, if no distinction between even and odd coordinates is necessary, we simply write $(\eta_i)$. 

A supermanifold of dimension $n|0$ is an ordinary differentiable manifold of dimension $n$, so at any time we may (and should) test the soundness of our theory by setting $m=0$.

One way of constructing examples of supermanifolds is the following: If $M$ is an ordinary differentiable manifold and $E\to M$ a vector bundle, then $(M,\Gamma (\Lambda E))$ is a supermanifold, where $\Gamma(\Lambda E)$ is the sheaf of sections of the exterior bundle $\Lambda E\to M$. Note that $\Gamma(\Lambda E)$ possesses a natural $\Z$-grading; in regarding it as the structural sheaf of a supermanifold, we retain only the induced $\Z_2$-grading. Although the Theorem of Batchelor \cite{Batchelor} asserts that every supermanifold over a differentiable manifold $M$ is isomorphic to one constructed in this way from some vector bundle over $M$, we prefer the definition given above since this isomorphism is non-canonical.

\subsection{Tangent Sheaf and Vector Fields}\label{Tangent}

For a super $\R$-algebra $A$, we give the endomorphisms ${\End\,} A$ of $A$ the structure of a super vector space via the natural grading $$({\End\,} A)_0=\{\varphi\in {\End\,} A\mid \varphi(A_0)\subset A_0, \varphi(A_1)\subset A_1\}, $$ $$ ({\End\,} A)_1=\{\varphi\in {\End\,}A\mid \varphi(A_0)\subset A_1, \varphi(A_1)\subset A_0\}.$$ Recall that a homogeneous element $\varphi\in {\End\,} A$ is a homogeneous derivation if 
\begin{equation}\varphi(ab)=\varphi (a)\cdot b+(-1)^{|\varphi||a|}a\cdot \varphi (b)\end{equation} for all $a,b\in A$, where for a homogeneous element $x$ of some graded object, $|x|\in \{0,1\}$ denotes the parity of $x$. An element $\varphi\in {\End\,} A$ is a derivation if its homogeneous components are homogeneous derivations.

If $M$ is a supermanifold of dimension $n|m$, we define  
$$\shT_M(U): ={\Der} (\shO_M(U)),$$ the $\shO_M(U)$-super module of derivations of $\shO_M(U)$. For $V\subset U\subset M$ there is a natural restriction map $\shT_M(U)\to \shT_M(V)$ turning $\shT_M$ into a sheaf of \mbox{$\shO_M$-super} modules, see \cite{Schmitt}, p.160. The $\shO_M$-module $\shT_M$ is locally free of dimension $n|m$, cf. \cite{DelMor}, $\S 3.3$. The sections of $\shT_M$ are called {\it vector fields}. We will refer to $\shT_M$ itself either as the {\it tangent sheaf} or the {\it tangent bundle} of the supermanifold $M$; this is not too much abuse of language, as is pointed out in \cite{DelMor}, $\S 3.4$
%%%%%%%%%%%%%% IM ARTIKEL WEG
%It would be the same to define $\shT_M$ by $$\shT_M(U):={\Der\, }(\shO_M|_U)=(\shDer\, \shO_M)(U)\subset (\shEnd\, \shO_M) (U) $$
%as a subsheaf of the sheaf of endomorphisms of $\shO_M$ (of course, these are endomorphisms of a sheaf of super vector spaces).
%%%%%%%%%%%%%% BIS HIER

On $\shT_M(U)$ we have a bracket $[\cdot,\cdot]$ defined by $$[X,Y]f:=X(Yf)-(-1)^{|X||Y|}Y(Xf).$$ It satisfies the {\it graded Jacobi identity}
 \begin{equation}[X,[Y,Z]]=[[X,Y],Z]+(-1)^{|X||Y|}[Y,[X,Z]]\label{Jacobi}\end{equation} and thus turns $\shT_M(U)$ into a Lie superalgebra.
 
For every point $p$ of $M$, the {\it tangent space} $T_pM$ of $M$ at $p$ is defined to be the space of derivations $\varphi:\shO_{M,p}\to \R$, i.e.~$$\varphi(fg)=\varphi(f) g(p)+(-1)^{|\varphi||f|}f(p)\varphi(g),$$where $\shO_{M,p}$ is the stalk of $\shO_M$ at $p$. For $p\in U$, there is a natural mapping $\shT_M(U)\to T_pM$ sending a vector field $X$ to its value $X_p$ at $p$; nevertheless, a vector field is not determined by its values at all points. The tangent spaces are the fibres of a bundle $TM\to M_{\red}$ of rank $n+m$ which canonically splits as the direct sum of the tangent bundle of the reduced manifold, $TM_{\red}\to M_{\red}$, and a bundle $(TM)_1\to M_{\red}$ with the odd parts of the tangent spaces as fibres. For a vector field $X$ on $M$, we denote by $\tilde{X}$ the associated section of $TM$.

The {\it cotangent bundle} of a supermanifold $M$ is by definition the dual $\Omega^1_M$ of $\shT_M$. As in \cite{DelMor}, we will write the duality pairing between the tangent and cotangent bundle as 
$$\left<\cdot,\cdot\right>:\shT_M\otimes_{\shO_M} \Omega^1_M\to \shO_M$$ with $\left<uX,v\omega\right>=(-1)^{|X||v|}uv\left< X,\omega\right>$ for $u,v\in \shO_M$. Then, as usual, we define $d:\shO_M\to \Omega^1_M$ by
$$\left< X,df\right> =Xf.$$

%%%%% Andere Definition d\Phi!!!!!

%Let $\Phi:M\to N$ be a morphism of supermanifolds. We want to show that it induces a morphism $d\Phi:\shT_M\to \phi^*\shT_N$ of $\shO_M$-modules (cf. \cite{Manin}, p. 187). \marginpar{why? Def. erkl\"aren} There is a well-defined pull-back morphism of $\shO_N$-modules \marginpar{pr\"ufen} $$\phi^*:\Omega^1_N\to \phi_*\Omega^1_M$$ given by $$df\mapsto d(\phi^*f);$$ thus, we may define $d\Phi$ as the composition
%$$d\Phi: \shT_M\to \shHom_{\shO_M}(\phi^*\Omega^1_N,\shO_M)=\phi^*\shT_N;\quad X\mapsto \{\alpha\mapsto \left<X,\phi^*\alpha\right>\}.$$

%Note that on the level of stalks, the morphism 
%\begin{equation}
%d\Phi_x:\shT_{M,x}\to (\phi^*\shT_N)_x=\shT_{N,\phi(x)}\otimes_{\shO_{N,\phi(x)}} \shO_{M,x}.
%\end{equation} is given by

A {\it vector field along a morphism} $\Phi:M\to N$ on $U\subset N$ is a morphism of super vector spaces $$X: \shO_N(U)\to \phi_*\shO_M(U)=\shO_M(\phi^{-1}(U))$$ such that its homogeneous components satisfy the derivation property \begin{equation} X(fg)=(Xf)\cdot \phi^*(g) + (-1)^{|X||f|}\phi^*(f)\cdot (Xg) \end{equation} for all $f,g\in \shO_N(U)$, cf. \cite{CariFig}. The set of such vector fields will be denoted by $\Der_\Phi (U)$ and the corresponding sheaf of vector fields along $\Phi$ by $\shT_\Phi:=\Der_\Phi$.

%As above, we could equivalently define a vector field along $\Phi$ as a local morphism of sheaves $\shO_N\to \phi_*\shO_M$, cf. \cite{CariFig1} \marginpar{genauer}

There are two standard ways of constructing vector fields along $\Phi$: If $X$ is a vector field on $N$, then $$\hat{X}:=\phi^*\circ X$$ is a vector field along $\Phi$; a vector field $Y$ on $M$ yields one by attaching $\phi^*$ on the other side: we define%\marginpar{d's vergleichen} 
\begin{equation} d\Phi(Y):=Y\circ \phi^*.
\end{equation}
If $\Phi$ is a diffeomorphism, we often use the same notation for the vector field $d\Phi(Y):=(\phi^{-1})^*\circ Y\circ \phi^*$ on $N$. 

If $\Phi:M\to N$ is a morphism of supermanifolds, then we have induced linear maps \mbox{$d_p\Phi:T_pM\to T_{\phi(p)}N$}. We call $\Phi$ an {\it immersion at} $p$ if $d_p\Phi$ is injective, and a {\it submersion at} $p$ if $d_p\Phi$ is surjective. See \cite{Varadara}, p.148 for the local structure of immersions and submersions.

The sheaf $\shT_\Phi$ is a locally free sheaf of $\phi_*\shO_M$-modules over $N$ of the same rank as $\shT_N$. More precisely, if $(x_i,\xi_\alpha)$ %\marginpar{Bez. Koord erkl\"aren: i, $\alpha$}
are local coordinates on $U\subset N$, then $(\phi^*\circ\partial_{x_i}=\hat{\partial}_{x_i}=\hat{\partial}_i,\phi^*\circ\partial_{\xi_\alpha}=\hat{\partial}_{\xi_\alpha}=\hat{\partial}_\alpha)$ is a basis of $\shT_\Phi(U)$: any $X\in \shT_\Phi(U)$ can be written uniquely as $$X=\sum f_i \hat{\partial}_{x_i}+\sum g_\alpha \hat{\partial}_{\xi_\alpha}$$ with $f_i,g_\alpha\in \phi_*\shO_M(U)$ \cite{CariFig}. Note that  
\begin{equation}\label{dPhiY in Koord} d\Phi(Y)=Y\circ \phi^*= \sum Y(\phi^*x_i)\cdot\hat{\partial}_{x_i}+\sum Y(\phi^*\xi_\alpha)\cdot\hat{\partial}_{\xi_\alpha}\end{equation} for all vector fields $Y$ on $M$.

\section{Lie Supergroups}
\subsection{Lie Supergroups and their Lie Superalgebras}

 %There are several ways to regard Lie supergroups -- we recall three of them. For the first two see e.g.~\cite{Varadara}.
 
A Lie supergroup is a group object in the category of supermanifolds, i.e.~a supermanifold $G$ together with morphisms $m:G\times G\to G,\, i:G\to G$ and $1:\R^{0|0}\to G$ representing the multiplication map, the inverse map and the unit element such that the usual group axioms are satisfied. The associativity law for example reads $$m\circ (\id_G\times m)=m\circ (m\times \id_G),$$ cf. Varadarajan \cite{Varadara}.

The reduced morphisms turn the reduced manifold $G_{\red}$ into a Lie group. 
A Lie supergroup $H$ is a {\it Lie subsupergroup} of a Lie supergroup $G$ if $H_{\red}$ is a Lie subgroup of $G_{\red}$ and the inclusion map of $H$ into $G$ is a morphism that is an immersion everywhere.

%Regarding the functor of points of a Lie supergroup we see that it takes values in the category of groups, the group structure on $S$-points being induced from the morphisms $m,i$ and $1$. Thus, we could equivalently define a Lie supergroup to be a representable contravariant functor from the category of supermanifolds to the category of groups.

% visualize Lie supergroups by means of Harish-Chandra pairs we have to define how to associate a Lie superalgebra to a given Lie supergroup. 

%\subsection{The Lie Superalgebra of a Lie Supergroup}
In classical Lie theory, a vector field $X$ on a Lie group $G$ is left-invariant if and only if $X_{gh}=dl_g(X_h)$ for all $g,h\in G$, where $l_g$ is left translation by $g$. We have to reformulate this in a way not using the elements of $G$ before generalizing the definition. A short calculation shows that it is equivalent to the condition \begin{equation}(I\otimes X)\circ m^*=m^*\circ X,\end{equation} where $I\otimes X$ is a vector field on $G\times G$, defined in the obvious way by acting only on the second component. This is now taken as the definition of {\it left-invariant vector field} on a Lie supergroup $G$; analogously, we say that $X$ is {\it right-invariant} if \begin{equation} (X\otimes I)\circ m^*=m^*\circ X.\end{equation} 
%%%%%%%%%%%%%% IM ARTIKEL WEG
%For any $g\in G_{\red}$, left translation with $g$ is a diffeomorphism $l_g:G\to G$ of the supergroup $G$. If $X$ is a left-invariant vector field on $G$, then 
%\begin{equation} X \circ l_g^*=l_g^*\circ X.\label{left-invariant_translation}\end{equation}
%%%%%%%%%%%%%% BIS HIER

The {\it Lie (super)algebra of $G$} is by definition the Lie superalgebra ${\mathfrak g}$ of all left-invariant vector fields on $G$. The usual isomorphism between the Lie algebra of $G$ and the tangent space of $G$ in the identity is still valid, as is proven in \cite{Varadara}, p.276f: The map ${\mathfrak g}\to T_eG;\, X\mapsto X_e$ is a linear isomorphism of super vector spaces. The converse map is also given: If $\tau\in T_eG$, \begin{equation}X_\tau:=(I\otimes \tau)\circ m^*\label{VectorToLInvVF}\end{equation} is the left invariant vector field with $(X_\tau)_e=\tau$; here, for a germ $f$ at $g\in G$, $m^*f$ is considered as a germ at $(g,e)$ so that $I\otimes \tau$ can be applied. 

We thus see that a Lie supergroup comes along with two objects we are more familiar with: its underlying Lie group and its Lie superalgebra. It would be nice if a Lie supergroup was already determined by this data.

\subsection{Harish-Chandra Pairs}\label{SectionHarishChandraPairs}

A {\it Harish-Chandra pair} is a pair $(G_0,\frak g)$, consisting of a Lie group $G_0$ and a Lie superalgebra $\frak g=\frak g_0\oplus\frak g_1$ with ${\rm Lie}(G_0)={\mathfrak g}_0$ and a representation $\Ad$ of $G_0$ on ${\mathfrak g}$ such that %\marginpar{Darst. durch Diff. in e eind. bestimmt? vgl. Koszul, Graded manifolds...}
\begin{enumerate}\item it extends the usual adjoint action of $G_0$ on its Lie algebra and
\item the differential of the action in the identity is equal to the Lie superbracket, restricted to ${\mathfrak g}_0\times {\mathfrak g}$. 
\end{enumerate}Note that if $G_0$ is connected, the first condition follows from the second.
%Note further that the second condition implies  \begin{equation}\Ad_{\exp(tX)}Y=X+t[X,Y]+\ldots\label{Ad_ad}\end{equation} for all $X\in {\mathfrak g_0}$ and $Y\in {\mathfrak g}$.

It is clear how to associate a Harish-Chandra pair to a given Lie supergroup $G$: Take the underlying Lie group $G_{\red}$, together with the Lie superalgebra $\frak g$ of $G$. The importance of the notion of Harish-Chandra pair results now from this functor being an equivalence of categories, cf. \cite{DelMor}, \S 3.8. 
%%%%%%%%%%%%% IM ARTIKEL WEG
%Let us explain why this is is a highly nontrivial statement: The odd part of the supermanifold structure of a Lie supergroup $G$ is in some sense trivial, see e.g.~\cite{Kostant}, Remark 3.6: \begin{equation}\shO_G(G_{\red})=\shC^\infty(G_{\red})\otimes \Lambda({\frak g}_1).\label{LieSupergroupGlobalTrivial}\end{equation} Nevertheless, it is a priori not obvious that the Harish-Chandra pair structure determines a (unique) Lie supergroup structure on $G$. One has to explicitly define the multiplication and inverse mappings; for these constructions, see \cite{Koszul} and \cite{BagSta}.
%%%%%%%%%%%%%% BIS HIER
See \cite{Koszul} and \cite{BagSta} for the construction of the Lie supergroup associated to a Harish-Chandra pair.

% this construction is described: Let $(G_0,{\mathfrak g})$ be a Harish-Chandra pair. For any open subset $U\subset G_0$ we define 
%\begin{equation}\shO_G(U):={\rm Hom}_{E({\mathfrak g}_0)}(E({\mathfrak g}),\shC^\infty_{G_0}(U)),\end{equation} where for a Lie (super)algebra $L$, $E(L)$ is the universal enveloping algebra of $L$, cf. \cite{Knapp} for the non-super case and \cite{Scheunert} for the super case. (Here, $E({\mathfrak g})$ and $\shC^\infty_{G_0}(U)$ are considered as $E({\mathfrak g}_0)$-modules.) Then $(G_0,\shO_G)$ is a supermanifold (more precisely, (\ref{LieSupergroupGlobalTrivial}) holds), roughly because $E({\mathfrak g})=E({\mathfrak g}_0)\otimes \Lambda({\mathfrak g}_1)$, cf. \cite{Koszul}, Lemma 1. Now in order to define the inverse morphism $i:G\to G$ we need to define for each open subset $U\subset G_0$ a morphism of super $\R$-algebras $i(U):\shO_G(U)\to \shO_G(U^{-1})$. This is done by  
%\begin{equation} (i^*(U)(\phi))(u)(g^{-1}):=\phi(\Ad_g \bar{u})(g)
%\end{equation} where $u\in E({\mathfrak g})$, $g\in U$ and $u\mapsto \bar{u}$ is the principal anti-automorphism $E({\mathfrak g})\to E({\mathfrak g})$, cf. \cite{Scheunert}, p. 21. 

%The multiplication morphism is given by \begin{equation} (m^*(U)(\phi))(u\otimes v)(g,h)=\phi(u\cdot \Ad_g v)(gh),\end{equation} see \cite{BagSta}, p. 583. \marginpar{Formeln verstehen}

%%%%%%%%%%%%%%% IM ARTIKEL WEG
%Note that if ${\mathfrak g}={\mathfrak g}_0\oplus {\mathfrak g}_1$ is a Lie superalgebra, there exists a unique simply-connected $G$ with Lie superalgebra ${\mathfrak g}$.

\subsection{Actions and Representations}\label{Section_Action}
An action of a Lie supergroup on a supermanifold is a morphism $G\times M\to M$ such that the usual axioms are satisfied.

In the language of Harish-Chandra pairs, such an action consists of an action of the reduced Lie group $G_{\red}$ on the supermanifold $M$, together with a morphism from $\mathfrak g$ to the opposite to the Lie superalgebra of vector fields on $M$, which are compatible in the sense that the differential of the action of the reduced group at the identity agrees with the restriction to the Lie algebra $\mathfrak g_0$, see \cite{DelMor}, p. 80.

If an action $\rho:G\times M\to M$ is given, this morphism ${\mathfrak g}\to \shT_M(M)^\circ$ is \begin{equation}X\mapsto (X_e\otimes I)\circ \rho^*;\label{ActionOddMorphism}\end{equation} note that it really is a morphism of Lie superalgebras, i.e.~%\marginpar{nachrechnen?}
\begin{equation}([X,Y]_e\otimes I)\circ \rho^*=-[(X_e\otimes I)\circ \rho^*,(Y_e\otimes I)\circ \rho^*].\label{BracketsSign}
\end{equation}

Similarly, a {\it representation} of a Lie supergroup $G$ on a super vector space $V$ consists of representations of $G_{\red}$ and ${\mathfrak g}$ on $V$ such that the differential of the representation of $G_{\red}$ coincides with the even part of the representation of ${\mathfrak g}$.

%An action $\psi$ of a Lie supergroup on a supermanifold is called {\it transitive} if...
%the induced actions on $S$-points are transitive. \marginpar{dr\"uber nachdenken...}

%Or: as in \cite{BoySan}: if there exists a point $p\in M$ such that the map $\psi_p:G\to M$ is an {\it epimorphism}, i.e.~the mapping $$\Hom(M,N)\to \Hom(G,N);\, \alpha\mapsto \alpha\circ \psi_p$$ is injective for all supermanifolds $N$.

\subsection{Examples}\label{section_LSG_examples}

In this section, we define those Lie supergroups and their Lie superalgebras that will be of importance for us in terms of their Harish-Chandra pairs. The adjoint action of the Lie group on the Lie superalgebra is always the standard one.

The {\it general linear supergroup} $\GL(n|m)$ is defined to be the Lie supergroup associated to the Harish-Chandra pair $$(\GL(n)\times \GL(m),\gl(n|m)),$$ where $\gl(n|m)$ is the Lie superalgebra consisting of block matrices $\left(\begin{array}{c|c}A & B \\ \hline C & D \end{array}\right)$ with $A,B,C$ and $D$ real $n\times n$-, $n\times m$-, $m\times n$- and $m\times m$-matrices, respectively. The gradation is given by $$\gl(n|m)_0=\{\left(\begin{array}{c|c}A & 0 \\ \hline 0 & D \end{array}\right)\}\text{ and }\gl(n|m)_1=\{\left(\begin{array}{c|c}0 & B \\ \hline C & 0 \end{array}\right)\}$$ and the bracket is the usual (anti-)commutator: for homogeneous elements $X,Y\in \gl(n|m)$ we define $[X,Y]:=XY-(-1)^{|X||Y|}YX$.

The {\it special linear supergroup} $\SL(n|m)$ is the Lie subsupergroup of $\GL(n|m)$ with reduced group $$\SL(n|m)_{\red}=\{(A,B)\in \GL(n)\times \GL(m)\mid {\det A}={\det B}>0\}$$ and Lie superalgebra $${\mathfrak{sl}}(n|m)=\{X\in \gl(n|m)\mid \str(X)=0\},$$ where the {\it supertrace} $\str$ of a matrix $\left(\begin{array}{c|c}A & B \\ \hline C & D \end{array}\right)$ is given by $\tr(A)-\tr(D)$. Note that the reduced group of $\SL(n|m)$ is isomorphic to $\SL(n)\times \SL(m)\times \R$ via $(A,B,\lambda)\mapsto (e^{\frac{\lambda}{n}}A, e^{\frac{\lambda}{m}}B)$.

In the case of $n=m$, the identity matrix $I_{2n}$ is an even element of the Lie superalgebra ${\mathfrak{sl}}(n|n)$ and generates a one-dimensional ideal; dividing this out, we obtain a Lie superalgebra denoted by $\psl(n|n)$. The corresponding Lie supergroup $\PSL(n|n)$ is given by the Harish-Chandra pair $$({\rm SL}(n|n)_{\red}/\R,\psl(n|n)).$$

The {\it orthosymplectic supergroup} $\OSp(n|2m)$ is the Lie subsupergroup of $\GL(n|2m)$ given by the Harish-Chandra pair $$(\rmO(n)\times \Sp(m;\R),\osp(n|2m)),$$ where
$$\osp(n|2m)=\{\left(\begin{array}{c|cc}A & B_1 & B_2 \\ \hline -B_2^t & C_1 & C_2 \\ B_1^t & C_3 & -C_1^t\end{array}\right)\mid A^t=-A,\, C_2^t=C_2,\, C_3^t=C_3\}.$$ Here, the real symplectic group $\Sp(m;\R)$ is the group of those transformations of $\R^{2m}$ leaving invariant the standard symplectic form. 
%%%%%%%%%%%%% IM ARTIKEL WEG
%The supergroup $\OSp(n|2m)$ should be thought of as the group of transformations of the super vector space $\R^{n|2m}$ leaving invariant the standard scalar superproduct, see \ref{Section_Metrics}. The name is justified by the fact that the reduced Lie group is the product of an orthogonal group and a symplectic group. 
The {\it special orthosymplectic supergroup} $\SOSp(n|2m)$ is the connected component of $\OSp(n|2m)$.

The {\it unitary supergroup} $\rmU(n|m)$ is the Lie supergroup associated to the Harish-Chandra pair
$$(\rmU(n)\times \rmU(m),{\mathfrak u}(n|m)),$$ where 
$${\mathfrak u}(n|m)=\{\left(\begin{array}{c|c}A & B \\ \hline -iB^* & C\end{array}\right)\mid A,B,C \text{ complex}, A^*=-A, C^*=-C\}.$$

\section{Riemannian Supergeometry}\label{chapter_RS}

\subsection{Metrics}\label{Section_Metrics}

A {\it scalar superproduct} on a super vector space $V=V_0\oplus V_1$ over a field $K$ (we will be interested almost exclusively in the case $K=\R$) is a non-degenerate graded-symmetric even $K$-bilinear form $\left< \cdot,\cdot\right>:V\times V\to K$. Here, the condition of graded symmetry is supposed to mean 
$$\left<X,Y\right>=(-1)^{|X||Y|}\left<Y,X\right>$$ for all homogeneous $X,Y\in V$. Since $K$ is considered as a purely even object, i.e.~$K_1=0$, being even means $\left<X,Y\right>=0$ for homogeneous $X,Y$ of different parity. Note that a graded scalar product on $V$ is the sum of a non-degenerate symmetric bilinear form $\left<\cdot,\cdot\right>_0$ on $V_0$ and a symplectic (i.e.~non-degenerate alternating bilinear) form $\left<\cdot,\cdot\right>_1$ on $V_1$. In particular, the existence of a scalar superproduct on $V$ forces $V_1$ to be even-dimensional.  
\begin{Rem}Note that a scalar superproduct on a purely even vector space is {\it not} the same as a scalar product since we do not impose any kind of positivity -- $\left<\cdot,\cdot\right>_0$ may be indefinite. See also \cite{Cortes}.
\end{Rem}
We define a {\it graded Riemannian metric} on a supermanifold $M$ as a graded-symmetric even non-degenerate $\shO_M$-linear morphism of sheaves 
$$\left<\cdot,\cdot\right>:\shT_M\otimes \shT_M\to \shO_M,$$ the non-degeneracy meaning that the mapping $X\mapsto \left<X,\cdot\right>$ is an isomorphism $\shT_M\to \Omega_M^1$. A supermanifold equipped with a graded Riemannian metric is called a {\it Riemannian supermanifold}.
 
For each $p\in M$ the morphism $\left<\cdot,\cdot\right>$ defines a scalar superproduct  $\left<\cdot,\cdot\right>_p$ on the real super vector space $T_p M$. As always, this family of scalar superproducts does not determine the graded Riemannian metric unless the odd dimension of $M$ is zero.
%A graded metric is called {\it graded Riemannian metric} if for every $p\in M$ the induced bilinear form $\left<\cdot,\cdot\right>_{p}$ on $T_pM$ is a scalar superproduct. 
%Since the positive linear combination of symplectic forms may become degenerate, the positive linear combination of graded Riemannian metrics need not be a metric. Hence, the usual proof of the existence of a Riemannian metric by partition of unity can not be applied. \marginpar{Example?}

Let us do the usual consistency check: A graded Riemannian metric induces in a natural way a pseudo-Riemannian metric on the underlying manifold $M_{\red}$: take $\left<\cdot,\cdot\right>_{p,0}$ on $T_pM_{\red}=(T_pM)_0$. Furthermore, on a usual differentiable manifold, the notion of graded Riemannian metric equals the notion of pseudo-Riemannian metric.

\begin{Rem} 
The choice of name for our metrics is justified by the fact that the natural metrics on Lie supergroups and homogeneous superspaces induced by the Killing form almost never have a definite sign, see \ref{SectionKilling}. The class of metrics on supermanifolds such that the underlying manifold is Riemannian (and not only pseudo-Riemannian) seems to play only a minor role. See also \cite{MontSan} and \cite{MontSan2}.

We also remark that there is another way of defining metrics having the nice property that the positive linear combination of metrics again is a metric -- which in our context obviously is not fulfilled. By passing to complexifications one may apply a definition of Tuynman in the context of cs manifolds, see \cite{Tuynman}, p. 188ff.
\end{Rem}

%%%%%%%%%%%% IM ARTIKEL WEG
%Locally, metrics can be represented in matrix form: if $(x_i,\xi_\alpha)$ are coordinates on $U\subset M$, the metric on $U$ is given by the matrix \begin{equation}
%\left(\begin{array}{c|c} (\left<\partial_i,\partial_j\right>)_{i,j} & (\left<\partial_i,\partial_\alpha\right>)_{i,\alpha} \\ \hline (\left<\partial_\alpha,\partial_i\right>)_{\alpha,i} & (\left<\partial_\alpha,\partial_\beta\right>)_{\alpha,\beta}\end{array}\right)=\left(\begin{array}{c|c} A & B \\ \hline B^t & C\end{array}\right)\label{MetricInCoord}
%\end{equation} with $A,B$ and $C$ matrices of superfunctions on $U$ satisfying $A=A^t$ and \mbox{$C=-C^t$}. 
%%%%%%%%%%%% BIS HIER

If a morphism $\Phi:M\to N$ and a graded Riemannian metric on $N$ are given, we can naturally evaluate vector fields along $\Phi$ with the metric via the morphism $\shT_\Phi\otimes \shT_\Phi\to \phi_*\shO_{M}$ given in coordinates $(\eta_i)$ by
\begin{equation}\left<\hat{\partial}_i,\hat{\partial}_j\right>=\left<\phi^*\circ \partial_i,\phi^*\circ \partial_j\right>:=\phi^*\left<\partial_i,\partial_j\right>.\label{MetricAlongPhi}
\end{equation} %\marginpar{well-defined?}

\subsection{Connections}

Let $(M,\shO_M)$ be a supermanifold and $\shE$ a locally free sheaf of $\shO_M$-super modules on $M$. A {\it connection} on $\shE$ (cf. \cite{DelMor}) is an even morphism \mbox{$\nabla:\shE\to \Omega^1_M\otimes \shE$} of sheaves of $\R$-super vector spaces that satisfies the Leibniz rule %\marginpar{geht hier schon die Parit\"at von $\nabla$ ein?}
\begin{equation}\nabla (fv)=df\otimes v+f\nabla v\end{equation} for all sections $f$ of $\shO_M$ and $v$ of $\shE$. If  we define \begin{equation}\nabla_Xv:=\left< X,\nabla v\right>\end{equation} for any vector field $X$, where $\left<X,\alpha\otimes v\right>:=\left<X,\alpha\right>v$, we get  
\begin{equation}\nabla_Xfv=Xf\cdot v+(-1)^{|X||f|}f\nabla_Xv \quad\text{and}\quad  |\nabla_X v|=|X|+|v|.\end{equation} 

In the case $\shE=\shT_M$ (in this case we speak of a connection on $M$) we define the {\it torsion} of a connection $\nabla$ on $\shT_M$ by \begin{equation} T_\nabla(X,Y):=\nabla_X Y-(-1)^{|X||Y|}\nabla_Y X -[X,Y].\end{equation} An easy calculation shows that for any superfunction $f$,  $$T_\nabla(fX,Y)=(-1)^{|f||X|}T_\nabla(X,fY)=fT_\nabla(X,Y).$$ Note that the tensorial property, namely that the values of the tensor depend only on the values of the inserted vector fields, is true in this context; nevertheless, it is less useful since a vector field cannot be reconstructed from its values.

The {\it curvature} of $\nabla$ is by definition \begin{equation}R(X,Y)Z:= [\nabla_X,\nabla_Y]Z-\nabla_{[X,Y]}Z.\end{equation} 

If $M$ is furnished with a graded Riemannian metric $\left<\cdot,\cdot\right>$, we call a connection $\nabla$ {\it metric} if 
\begin{equation}X\left<Y,Z\right>=\left<\nabla_X Y,Z\right>+(-1)^{|X||Y|}\left<Y,\nabla_X Z\right>.\label{eqMetricConnection}\end{equation}
We also have the usual notion of {\it Christoffel symbols}: If $(\eta_i)$ is a system of coordinates (both even and odd) on $U\subset M$,  
\begin{equation} \nabla_{\partial_i}\partial_j=\sum_k \Gamma_{ij}^k \partial_k
\end{equation} gives well-defined elements $\Gamma_{ij}^k\in \shO_M(U)$ of parity \begin{equation}\label{ParityChristoffel} |\Gamma_{ij}^k|=|\eta_i|+|\eta_j|+|\eta_k|.
\end{equation} If it is necessary to distinguish between the even and the odd coordinates, in order not to let the notation explode the indices (latin or greek) have to indicate which Christoffel symbol is meant, e.g.~\mbox{$\nabla_{\partial_{\xi_\alpha}} \partial_{x_i}=\sum_{j}\Gamma_{\alpha i}^j \partial_{x_j}+\sum_{\beta} \Gamma_{\alpha i}^\beta \partial_{\xi_\beta}$}.

Note that a connection on $M$ in a natural way induces a connection on the vector bundle $TM\to M_{\red}$ by reduction of the coefficient functions.

The following theorem, whose proof is -- apart from the additional signs -- the same as in standard (pseudo-)Riemannian geometry, can be found in \cite{MontSan}: 
\begin{Thm} On a supermanifold $M$ with a graded Riemannian metric, there exists a unique torsionless and metric connection $\nabla$ (which will be called the {\emph Levi-Civita connection} of the metric). It is implicitly defined by the formula 
\begin{align} 2\left<\nabla_X Y,Z\right> &= X\left<Y,Z\right> -(-1)^{|Z|(|X|+|Y|)}Z\left< X,Y\right> 
\nonumber \\ 
 &\quad\quad +(-1)^{|X|(|Y|+|Z|)}Y\left< Z,X\right>+\left< [X,Y],Z\right> \label{LeviCivita_Formula} \\ &\qquad-(-1)^{|X|(|Y|+|Z|)} \left< [Y,Z],X\right> \nonumber \\
 &\qquad+  (-1)^{|Z|(|X|+|Y|)}\left<[Z,X],Y\right>. \nonumber
\end{align}
\end{Thm}

The Levi-Civita connection of the metric induces the standard Levi-Civita connection of the induced pseudo-Riemannian metric on $M_{\red}$, as can be seen by regarding (\ref{LeviCivita_Formula}) on the level of the even tangent spaces $(T_pM)_0$ -- there, the formula becomes the usual formula defining the Levi-Civita connection.

Consequently, if we are in coordinates $(x_i,\xi_\alpha)$ on $U$, the projections of the Christoffel symbols $\Gamma_{ij}^k$ on $\shC^\infty(U)$ are the Christoffel symbols of the Levi-Civita connection on $M_{\red}$ with respect to the corresponding coordinates on $M_{\red}$.

\begin{Prop} \label{CurvProperties}Let $M$ be a supermanifold with a graded Riemannian me\-tric and the corresponding Levi-Civita connection. Then the following equalities hold for all vector fields $X,Y,Z,W$:
\begin{equation}\left<R(X,Y)Z,W\right>=-(-1)^{|X||Y|}\left<R(Y,X)Z,W\right>=-(-1)^{|Z||W|}\left<R(X,Y)W,Z\right>\label{CurvSymmetry1}\end{equation}
\begin{equation}\left<R(X,Y)Z,W\right>=(-1)^{(|X|+|Y|)(|Z|+|W|)}\left<R(Z,W)X,Y\right>\label{CurvSymmetry2}
\end{equation}
\begin{equation} R(X,Y)Z+(-1)^{|Z|(|X|+|Y|)}R(Z,X)Y+(-1)^{|X|(|Y|+|Z|)}R(Y,Z)X=0.\label{1stBianchi}
\end{equation}
\end{Prop}
\begin{proof} Except for the additional signs, these relations are proven just like in the standard theory; as an example, we calculate one part of (\ref{CurvSymmetry1}).
\begin{align*}\left<R(X,Y)Z,W\right>&=\left<\nabla_X\nabla_YZ,W\right>-(-1)^{|X||Y|}\left<\nabla_Y\nabla_XZ,W\right>-\left<\nabla_{[X,Y]}Z,W\right>\\
&=X\left<\nabla_YZ,W\right>-(-1)^{|X|(|Y|+|Z|)}\left<\nabla_YZ,\nabla_XW\right>\\
&\quad\quad-(-1)^{|X||Y|}Y\left<\nabla_XZ,W\right>+(-1)^{|Y||Z|}\left<\nabla_XZ,\nabla_YW\right>\\
&\quad\quad-[X,Y]\left<Z,W\right>+(-1)^{|Z|(|X|+|Y|)}\left<Z,\nabla_{[X,Y]}W\right>\\
&=XY\left<Z,W\right>-(-1)^{|Y||Z|}X\left<Z,\nabla_YW\right>\\
&\quad\quad-(-1)^{|X|(|Y|+|Z|)}\left<\nabla_YZ,\nabla_XW\right>-(-1)^{|X||Y|}YX\left<Z,W\right>\\
&\quad\quad+(-1)^{|X|(|Y|+|Z|)}Y\left<Z,\nabla_XW\right>+(-1)^{|Y||Z|}\left<\nabla_XZ,\nabla_YW\right>\\
&\quad\quad-[X,Y]\left<Z,W\right>+(-1)^{|Z|(|X|+|Y|)}\left<Z,\nabla_{[X,Y]}W\right>\\
&=-(-1)^{|Z||W|}\left<R(X,Y)W,Z\right>.
\end{align*}
%\end{proof}
%Since $R$ is a tensor, it suffices to check the first Bianchi identity (\ref{1stBianchi}) for vector fields $X,Y,Z$ with $[X,Y]=[Y,Z]=[X,Z]=0$; we omit this calculation because it is straightforward.
%\begin{align*}&R(X,Y)Z+(-1)^{|Z|(|X|+|Y|)}R(Z,X)Y+(-1)^{|X|(|Y|+|Z|)}R(Y,Z)X\\
%\quad &\quad= \nabla_X\nabla_YZ-(-1)^{|X||Y|}\nabla_Y\nabla_XZ +(-1)^{|Z|(|X|+|Y|)}\nabla_Z\nabla_XY\\
%&\quad\quad\quad-(-1)^{|Z||Y|}\nabla_X\nabla_ZY+(-1)^{|X|(|Y|+|Z|)}\nabla_Y\nabla_ZX\\
%&\quad\quad\quad-(-1)^{|X||Y|+|X||Z|+|Y||Z|}\nabla_Z\nabla_YX\\
%&\quad=\nabla_X[Y,Z]-(-1)^{|X||Y|}\nabla_Y[X,Z]+(-1)^{|Z|(|X|+|Y|)}\nabla_Z[X,Y]=0.
%\end{align*} 
Note that (\ref{CurvSymmetry2}) is a consequence of (\ref{CurvSymmetry1}) and (\ref{1stBianchi}) -- this calculation is explicitly carried out in \cite{Cortes}, Lemma 1.
%%%%%%%%%% IM ARTIKEL WEG
%\begin{align*}\left<R(X,Y)Z,W\right>&=-(-1)^{|Z|(|X|+|Y|)}\left<R(Z,X)Y,W\right>\\
%&\quad\quad-(-1)^{|X|(|Y|+|Z|)}\left<R(Y,Z)X,W\right>\\
%&=(-1)^{|Z|(|X|+|Y|)+|Y||W|}\left<R(Z,X)W,Y\right>\\
%&\quad\quad+(-1)^{|X|(|Y|+|Z|+|W|)}\left<R(Y,Z)W,X\right>\\
%&=-(-1)^{(|X|+|Y|)(|Z|+|W|)+|W||Z|}\left<R(W,Z)X,Y\right>\\
%&\quad\quad-(-1)^{|Z||W|+|Z||Y|+|Y||W|}\left<R(X,W)Z,Y\right>\\
%&\quad\quad-(-1)^{(|X|+|W|)(|Y|+|Z|)+|W||X|}\left<R(W,Y)Z,X\right>\\
%&\quad\quad-(-1)^{(|X|+|Y|)(|Z|+|W|)+|X||Y|}\left<R(Z,W)Y,X\right>\\
%&=(-1)^{(|X|+|Y|)(|Z|+|W|)}2\left<R(Z,W)X,Y\right>\\
%&\quad\quad+(-1)^{|W|(|Y|+|Z|)}\left<R(X,W)Y,Z\right>\\
%&\quad\quad+(-1)^{|W|(|Y|+|Z|)+|X|(|W|+|Y|)}\left<R(W,Y)X,Z\right>\\
%&=(-1)^{(|X|+|Y|)(|Z|+|W|)}2\left<R(Z,W)X,Y\right>\\
%&\quad\quad-(-1)^{|X||Y|+|W||Z|}\left<R(Y,X)W,Z\right>\\
%&=(-1)^{(|X|+|Y|)(|Z|+|W|)}2\left<R(Z,W)X,Y\right>-\left<R(X,Y)Z,W\right>
%\end{align*}
%%%%%%%%%%%% BIS HIER
\end{proof}

\subsection{Covariant Derivatives along Supercurves}

Let $\gamma=(\gamma,\gamma^*):\R^{1|1}\to (M,\shO_M)$ be a supercurve in a supermanifold with a connection $\nabla$. Let $(t,\xi)$ be the standard coordinates on $\R^{1|1}$. Recall (\ref{dPhiY in Koord}) that if $(\eta_i)$ are coordinates on $U\subset M$ (for the moment we do not distinguish between even and odd coordinates), then %\marginpar{Koord auf $\R^{1|1}$} 
$$d\gamma(\partial_t)=\partial_t\circ \gamma^*= \sum \partial_t(\gamma^*\eta_i)\cdot\hat{\partial}_{\eta_i}.$$
We define the covariant derivative of a vector field $X=\sum f_j \hat{\partial}_{\eta_j}$ along $\gamma$ with respect to the even coordinate by
\begin{align} \frac{\nabla}{dt}X&=\sum_j (\partial_t f_j)\cdot \hat{\partial}_{\eta_j}+f_j\frac{\nabla}{dt}\hat{\partial}_{\eta_j}\nonumber\\
&=\sum_j (\partial_t f_j)\cdot \hat{\partial}_{\eta_j}+f_j\sum_i \partial_t (\gamma^*\eta_i)\cdot \gamma^*\nabla_{\partial_{\eta_i}} \partial_{\eta_j}\nonumber\\
%&=\sum_j (\partial_t f_j)\cdot \hat{\partial}_{\eta_j}+f_j\sum_i \partial_t (\gamma^*\eta_i)\cdot \sum_k \gamma^*\Gamma_{ij}^k \cdot \hat{\partial}_{\eta_k}\nonumber\\
&=\sum_k \left( \partial_t f_k+\sum_{i,j}f_j\cdot \partial_t(\gamma^*\eta_i)\cdot \gamma^*\Gamma_{ij}^k\right)\hat{\partial}_{\eta_k}.\label{CovDerivativeCurveEven}
\end{align}
Note that the $\gamma^*\Gamma_{ij}^k$ are sections of $\shO_{\R^{1|1}}$ over $U$, i.e.~elements of $\shC^\infty(\gamma^{-1}(U))[\xi]$; furthermore, they are homogeneous and hence either elements of $\shC^\infty(\gamma^{-1}(U))$ or $\shC^\infty(\gamma^{-1}(U))\cdot \xi$, depending on the parity given by (\ref{ParityChristoffel}).

Analogously to (\ref{CovDerivativeCurveEven}), we can define the covariant derivative along a curve with respect to the odd coordinate -- we only have to take into account that it has to become an odd operator:
\begin{align} \frac{\nabla}{d\xi}X%&=\sum_j (\partial_\xi f_j)\cdot \hat{\partial}_{\eta_j}+(-1)^{|f_j|}f_j\frac{\nabla}{d\xi}\hat{\partial}_{\eta_j}\nonumber\\
%&=\sum_j (\partial_\xi f_j)\cdot \hat{\partial}_{\eta_j}+(-1)^{|f_j|}f_j\sum_i \partial_\xi (\gamma^*\eta_i)\cdot \gamma^*\nabla_{\partial_{\eta_i}} \partial_{\eta_j}\nonumber\\
%&=\sum_j (\partial_t f_j)\cdot \hat{\partial}_{\eta_j}+f_j\sum_i \partial_t (\gamma^*\eta_i)\cdot \sum_k \gamma^*\Gamma_{ij}^k \cdot \hat{\partial}_{\eta_k}\nonumber\\
&=\sum_k \left( \partial_\xi f_k+\sum_{i,j}(-1)^{|f_j|}f_j\cdot \partial_\xi(\gamma^*\eta_i)\cdot \gamma^*\Gamma_{ij}^k\right)\hat{\partial}_{\eta_k}.\label{CovDerivativeCurveOdd}
\end{align}

If $M$ is now assumed to be Riemannian and equipped with a metric connection, we get the following:
\begin{Prop} For vector fields $X$ and $Y$ along a supercurve $\gamma$, we have
$$\partial_t\left<X,Y\right>=\left<\frac{\nabla}{dt}X,Y\right>+\left<X,\frac{\nabla}{dt}Y\right>\quad \text{ and}$$
$$\partial_\xi\left<X,Y\right>=\left<\frac{\nabla}{d\xi}X,Y\right>+(-1)^{|X|}\left<X,\frac{\nabla}{d\xi}Y\right>.$$
\end{Prop}
\begin{proof} Formula (\ref{dPhiY in Koord}) and the metricity yield \begin{align*} \partial_t\left<\hat{\partial}_i,\hat{\partial}_j\right>&=\partial_t\circ \gamma^*\left<\partial_i,\partial_j\right>=\sum_k\partial_t(\gamma^*\eta_k)\cdot\gamma^*\circ \partial_k\left<\partial_i,\partial_j\right>\\
&=\sum_k \partial_t(\gamma^*\eta_k)\cdot \gamma^*(\left<\nabla_{\partial_k}\partial_i,\partial_j\right>+(-1)^{|\eta_i||\eta_k|}\left<\partial_i,\nabla_{\partial_k}\partial_j\right>)\\
&=\left<\frac{\nabla}{dt}\hat{\partial}_i,\hat{\partial_j}\right>+\left<\hat{\partial}_i,\frac{\nabla}{dt}\hat{\partial}_j\right>;
\end{align*} the second equality is proven analogously.
%%%%%%% IM ARTIKEL WEG
%\begin{align*} \partial_\xi\left<\hat{\partial}_i,\hat{\partial}_j\right>&=\partial_\xi\circ \gamma^*\left<\partial_i,\partial_j\right>=\sum_k\partial_\xi(\gamma^*\eta_k)\cdot\gamma^*\circ \partial_k\left<\partial_i,\partial_j\right>\\
%&=\sum_k \partial_\xi(\gamma^*\eta_k)\cdot \gamma^*(\left<\nabla_{\partial_k}\partial_i,\partial_j\right>+(-1)^{|\eta_i||\eta_k|}\left<\partial_i,\nabla_{\partial_k}\partial_j\right>)\\
%&=\left<\frac{\nabla}{d\xi}\hat{\partial}_i,\hat{\partial_j}\right>+(-1)^{|\partial_i|}\left<\hat{\partial}_i,\frac{\nabla}{d\xi}\hat{\partial}_j\right>.
%\end{align*}
%%%%%%%% BIS HIER
\end{proof}

\subsection{Geodesics}

Before giving with (\ref{GeodesicsOK})
the right condition for being a supergeodesic we start with some false ones at first sight resembling the known condition from Riemannian geometry. 

Let a supermanifold $M$ be equipped with a graded Riemannian metric and the corresponding Levi-Civita connection. One generalization of the ordinary definition of geodesics would be to demand the vanishing of the term $\frac{\nabla}{dt}(d\gamma(\partial_t))$ and possibly additionally of some of the three terms  $\frac{\nabla}{dt}(d\gamma(\partial_\xi))$, $\frac{\nabla}{d\xi}(d\gamma(\partial_t))$ and $\frac{\nabla}{d\xi}(d\gamma(\partial_\xi))$. Let us first regard supercurves $\gamma$ that satisfy only $\frac{\nabla}{dt}(d\gamma(\partial_t))=0$. In local coordinates, this condition is
\begin{equation} \partial_t^2 (\gamma^*\eta_k)+\sum_{i,j}\partial_t(\gamma^*\eta_j)\cdot \partial_t(\gamma^*\eta_i)\cdot \gamma^*\Gamma_{ij}^k=0
\end{equation} for all $k$. 
%(As far as I know the only references for the definition of supergeodesics are \cite{DeWitt} and \cite{MonMun}.) 
It is quite useful to write these equations separately for the even and the odd coordinates. If such $(x_i,\xi_\alpha)$ are fixed, we have $\gamma^*x_i=g_i$ and $\gamma^*\xi_\alpha=h_\alpha\xi$ for some $g_i,h_\alpha\in \shC^\infty(\gamma^{-1}(U))$. Then for the even coordinates, the geodesic equations are
\begin{align} 0&= g_k''+\sum_{i,j}g_i'g_j' \gamma^*\Gamma_{ij}^k+\sum_{i,\beta}g_i'h_\beta'\underbrace{\xi\gamma^*\Gamma_{i\beta}^k}_{=0} \nonumber\\
&\quad\quad\,\,+\sum_{\alpha,j} h_\alpha' g_j'\underbrace{\xi\gamma^*\Gamma_{\alpha j}^k}_{=0}+\sum_{\alpha,\beta}h_\alpha'h_\beta'\underbrace{\xi^2}_{=0}\gamma^*\Gamma_{\alpha\beta}^k\nonumber\\
& =g_k''+\sum_{i,j}g_i'g_j'\gamma^*\Gamma_{ij}^k, \label{EvenGeodesicEquation}
\end{align} where we used (\ref{ParityChristoffel}) for the second and third sum -- there, we have Christoffel symbols with an odd number of greek indices. These are the usual geodesic equations of the underlying Riemannian manifold, which is an argument in favor of the definition: a geodesic in a supermanifold should be a supercurve with the underlying curve being an ordinary geodesic and satisfying some kind of additional odd condition. 

We also write down the equations for the odd coordinates:
\begin{align}0%&=h_\delta''\xi+\sum_{i,j}g_i'g_j' \gamma^*\Gamma_{ij}^\delta+\sum_{i,\beta}g_i'h_\beta'\xi\gamma^*\Gamma_{i\beta}^\delta+\sum_{\alpha,j} h_\alpha' g_j'\xi\gamma^*\Gamma_{\alpha j}^\delta \nonumber \\
&= h_\delta''\xi+\sum_{i,j}g_i'g_j' \gamma^*\Gamma_{ij}^\delta+2\sum_{i,\beta}g_i'h_\beta'\xi\gamma^*\Gamma_{i\beta}^\delta.  \label{OddGeodesicEquation} \end{align} In each of the summands appears exactly one $\xi$ so we reduced the equations to ordinary differential equations of second order.

What are the initial conditions for geodesics of this kind in coordinate-free notation? The value of $\gamma:\R^{1|1}\to M$ at $0$ is (in coordinates) given by the tuple $(g_i(0))$. Let us write $d_0\gamma:T_0 \R^{1|1}\to T_{\gamma(0)}M$ in coordinates: \begin{align}d\gamma(\left.\partial_t\right|_{0})=\left.\partial_t\right|_{0}\circ \gamma^*=\sum_{i}g_i'(0) \left.\partial_{x_i}\right|_{\gamma(0)} \label{Differential_Curve_Even}\\
d\gamma(\left.\partial_{\xi}\right|_{0})=\left.\partial_{\xi}\right|_{0}\circ \gamma^*=\sum_{\alpha}h_{\alpha}(0)\left.\partial_{\xi_\alpha}\right|_{\gamma(0)}
\label{Differential_Curve_Odd}\end{align}
where we used that the values of nilpotent functions are zero at any point. We thus see that such a geodesic is not determined by its value and the value of its differential at one point. One might argue that this is unavoidable because of the unimportant role of points in the world of supermanifolds, but then one would have to accept that the initial conditions can not be written in a natural way in coordinate-free notation. It would thus be favourable if the odd geodesic equations (\ref{OddGeodesicEquation}) were differential equations of first order. 
%%%%%%%%%%%%% IM ARTIKEL WEG
%Before we can achieve this by taking into account the bundle structure of supermanifolds, we have to make some calculations analogous to (\ref{EvenGeodesicEquation}) and (\ref{OddGeodesicEquation}):
%\begin{align}
%\frac{\nabla}{dt} d\gamma(\partial_\xi)&=\sum_k \left(\sum_{i,\beta}h_\beta g_i'\gamma^*\Gamma_{i\beta}^k+\sum_{\alpha,\beta}h_\beta h_\alpha'\xi\gamma^*\Gamma_{\alpha\beta}^k\right)\hat{\partial}_{x_k}\nonumber\\
%&\qquad +\sum_\delta \left(h_\delta'+\sum_{i,\beta} h_\beta g_i'\gamma^*\Gamma_{i\beta}^\delta\right)\hat{\partial}_{\xi_\delta}.\label{CovCurve_t_xi}\\
%\frac{\nabla}{d\xi} d\gamma(\partial_t)&=\sum_k \left(\sum_{\alpha,j}g_j'h_\alpha \gamma^*\Gamma_{\alpha j}^k-\sum_{\alpha,\beta}h_\beta h_\alpha'\xi\gamma^*\Gamma_{\alpha\beta}^k\right)\hat{\partial}_{x_k}\nonumber\\
%&\qquad +\sum_\delta \left(h_\delta'+\sum_{\alpha,j} g_j'h_\alpha \gamma^*\Gamma_{\alpha j}^\delta\right)\hat{\partial}_{\xi_\delta}.\label{CovCurve_xi_t}\\
%\frac{\nabla}{d\xi} d\gamma(\partial_\xi)&=0.\label{CovCurve_xi_xi}
%\end{align}

The additional vanishing of some of the three terms $\frac{\nabla}{dt}d\gamma(\partial_\xi)$, $\frac{\nabla}{d\xi}d\gamma(\partial_t)$ and $\frac{\nabla}{d\xi}d\gamma(\partial_\xi)$ does not solve the problem described above, whereas it is solved via the following definition - recall that the tilde means passing to the reduced level.
\begin{Def}\label{Def_SG}A curve $\gamma:\R^{1|1}\to M$ is a {\emph {(super)geodesic}} if
\begin{equation}\left(\frac{\nabla}{dt}d\gamma(\partial_t)\right)^{\widetilde{}}=\left(\frac{\nabla}{dt} d\gamma(\partial_\xi)\right)^{\widetilde{}}=0.\label{GeodesicsOK}\end{equation} 
\end{Def}
Easy calculations in coordinates show that for all supercurves $\gamma$, we have $\frac{\nabla}{d\xi}d\gamma(\partial_\xi)=0$ and $\big(\frac{\nabla}{dt}d\gamma(\partial_\xi)\big)^{\widetilde{}}=\big(\frac{\nabla}{d\xi}d\gamma(\partial_t)\big)^{\widetilde{}}.$

\begin{Prop} \label{GeodesicDetermined}A curve $\gamma:\R^{1|1}\to M$ is a supergeodesic if and only if the underlying curve $\widetilde{\gamma}:\R\to M$ is a geodesic and $\left(\frac{\nabla}{dt}d\gamma(\partial_\xi)\right)^{\widetilde{}}=0$; in coordinates, the second condition is  \begin{equation}h_\delta'+\sum_{i,\beta} g_i'h_\beta \gamma^*\Gamma_{i\beta}^\delta=0\label{RightOddGeodesicEquation}\end{equation} for all $\delta$ (the so-called {\emph {odd geodesic equations}}).

Furthermore, for every $p\in M$ and every tangent vector $\tau\in T_pM$, there exists a unique supergeodesic $\gamma$ with $\gamma(0)=p$ and $d_0\gamma(\partial_t+\partial_\xi)=\tau$.\label{GeodesicConditions}
\end{Prop} 
\begin{proof}The first statement is clear by simply writing it down in coordinates. Since the odd differential equations (\ref{RightOddGeodesicEquation}) are of first order, the second statement follows by (\ref{Differential_Curve_Even}) and (\ref{Differential_Curve_Odd}).
\end{proof}

The formal consequences of this definition are thus just what we wanted them to be; on the other hand, this notion of supergeodesic seems highly unnatural compared to the ones suggested before. Nevertheless, it will turn out that it reflects the bundle structure of a supermanifold in a very natural way.

If a supermanifold $M$ is given, the Theorem of Batchelor allows us to find a vector bundle $\pi:E\to M_{\red}$ such that $M$ is given by the sheaf of sections of the exterior bundle $\Lambda E\to M_{\red}$, cf. \ref{Section_Smfd}.
The data of a supercurve \mbox{$\gamma=(\widetilde{\gamma},\gamma^*):\R^{1|1}\to M$} is now the same as the data of an ordinary curve $\bar{\gamma}:\R\to E$: If a local basis $\xi_\alpha$ of $E$ is given (which at the same time is a local system of odd coordinates of $M$), then $\gamma^*\xi_\alpha=g_\alpha\xi$ for some $\shC^\infty$-functions $g_\alpha$. Then we define a lift $\bar{\gamma}(t):\R\to E$ of $\widetilde{\gamma}$ by $\bar{\gamma}(t)=\sum_\alpha g_\alpha(t){\xi_\alpha}({\widetilde{\gamma}(t)})$; note that the expression $\xi_\alpha(\widetilde{\gamma}(t))$ makes sense since we regard the $\xi_\alpha$ as local sections of $E$. Clearly, this correspondence between supercurves in $M$ and curves in $E$ is one-to-one.

We want to justify Definition \ref{Def_SG} by proving that the property of $\gamma$ being a supergeodesic is equivalent to a natural property of the associated curve $\bar{\gamma}$. The vector bundle $E$ carries a natural connection induced from the connection on the supermanifold $M$, defined as follows: Let $(x_i,\xi_\alpha)$ be coordinates on $M$ such that the $\xi_\alpha$ are a local basis of the vector bundle $E$. For a vector field $X$ on $M_{\red}$, we define $$\bar{\nabla}_X \xi_\alpha:=\left(\nabla_X \partial_{\xi_\alpha}\right)^E,$$ where the superscript $E$ is supposed to mean the following: write the odd vector field $\nabla_X\partial_{\xi_\alpha}$ in coordinates, reduce the coefficient functions and replace $\partial_{\xi_\beta}$ by $\xi_\beta$ to arrive at a section of $E$. As an example, $(\xi_1\partial_x+(x^2+\xi_1\xi_2)\partial_{\xi_2})^E=x^2\xi_2$. Note that $\bar{\nabla}$ really is a connection on $E$. Now Proposition \ref{GeodesicConditions} shows that $\gamma$ is a supergeodesic if and only if $\bar{\gamma}$ is the horizontal lift of a geodesic in $M_{\red}$, where horizontality is meant with respect to $\bar{\nabla}$. 

\subsection{Parallel Displacement}
A vector field $X$ along a supercurve $\gamma$ in a Riemannian supermanifold $M$, equipped with the Levi-Civita connection, is called {\it parallel} if \begin{equation}\left(\frac{\nabla}{dt}X\right)^{\widetilde{}}=\left(\frac{\nabla}{d\xi}X\right)^{\widetilde{}}=0.\label{VFieldParallel}\end{equation} 
\begin{Prop} For each tangent vector $\tau\in T_{\gamma(0)}M$ there is a unique pa\-rallel vector field $X$ along $\gamma$ with value $\tau$ at $0$. Furthermore, if $\tau$ is homogeneous, $X$ is homogeneous of the same parity.
\end{Prop}
\begin{proof} Writing the conditions (\ref{VFieldParallel}) for a vector field $$X=\sum_k (f_j+\xi g_j)\hat{\partial}_j+\sum_\beta (f_\beta+\xi g_\beta)\hat{\partial}_\beta$$ along $\gamma$ (where $f_j,g_j,f_\beta$ and $g_\beta$ are smooth functions) in coordinates $(x_i,\xi_\alpha)$, we see that $X$ is parallel if and only if 
\begin{align*} f_k'+\sum_{i,j} f_j\cdot (\partial_t \gamma^*x_i)^{\widetilde{}}\cdot(\gamma^*\Gamma_{ij}^k)^{\widetilde{}}=g_k+\sum_{\alpha,\beta}f_\beta\cdot(\partial_\xi \gamma^*\xi_\alpha)^{\widetilde{}}\cdot(\gamma^*\Gamma_{\alpha\beta}^k)^{\widetilde{}}=0
\end{align*} for all $k$ and 
\begin{align*} f_\delta'+\sum_{i,\beta}f_\beta\cdot (\partial_t \gamma^*x_i)^{\widetilde{}}\cdot (\gamma^*\Gamma_{i\beta}^\gamma)^{\widetilde{}}=g_\delta+\sum_{\alpha,j} f_j\cdot (\partial_\xi\gamma^*\xi_\alpha)^{\widetilde{}}\cdot(\gamma^*\Gamma_{\alpha j}^\gamma)^{\widetilde{}}=0
\end{align*} for all $\delta$. The initial condition is given by the tuples $(f_j(0))$ and $(f_\beta(0))$, so the unique existence follows. If $\tau$ is even, all the $f_\beta(0)$ are equal to $0$. Thus, the $f_\beta$ and $g_j$ vanish completely and hence $X$ is even. If $\tau$ is odd, all the $f_j(0)$ are equal to $0$ and thus the $f_j$ and $g_\beta$ vanish, so $X$ is odd.
\end{proof}
We thus may define {\it parallel displacement} $P(\gamma)_s^t:T_{\gamma(s)}M\to T_{\gamma(t)}M$ just like in the standard theory: if $\tau\in T_{\gamma(s)}M$ is given, let $X$ be the unique parallel vector field along $\gamma$ with $X_s=\tau$ and define $P(\gamma)_s^t(\tau):=X_t$.
\begin{Prop} Parallel displacement is an isometry.
\end{Prop}
\begin{proof} For any two parallel vector fields $X$ and $Y$ along $\gamma$, the superfunction $$\partial_t\left<X,Y\right>=\left<\frac{\nabla}{dt}X,Y\right>+\left<X,\frac{\nabla}{dt}Y\right>$$ is nilpotent. Thus, the function $t\mapsto \left<X_t,Y_t\right>_{\tilde{\gamma}(t)}$ is constant.
\end{proof}

\subsection{Isometries} \label{SectionIsometries}

Let $M$ and $N$ be Riemannian supermanifolds, equipped with their Levi-Civita connections. We say that a diffeomorphism $\Phi=(\phi,\phi^*):M\to N$ is an {\it isometry} if it respects the metric: $\Phi^*\left<\cdot,\cdot\right>=\left<\cdot,\cdot\right>$, i.e.~$$\phi^*\left<d\Phi(X),d\Phi(Y)\right>=\phi^*\left<(\phi^{-1})^*\circ X\circ \phi^*,(\phi^{-1})^*\circ Y\circ \phi^*\right>=\left<X,Y\right>$$ for all vector fields $X,Y$ on $M$. (Recall the different definitions for $d\Phi(X)$ given in \ref{Tangent} -- since $\Phi$ is a diffeomorphism we regard it as a vector field on $N$.)  

If $\Phi:M\to N$ is an isometry, $\Phi$ is in particular {\it affine}, i.e.~$$d\Phi(\nabla_XY)=\nabla_{d\Phi(X)}d\Phi(Y)$$ for all vector fields $X,Y$ on $M$, as can be seen by regarding the formula defining the Levi-Civita connection, (\ref{LeviCivita_Formula}). Consequently, if $\eta_i$ are coordinates on $U$, $\Gamma_{ij}^k$ the Christoffel symbols with respect to these coordinates and $\bar{\Gamma}^k_{ij}$ the Christoffel symbols with respect to the coordinates $(\phi^{-1})^* \eta_i$ on $\phi(U)$, then ${\phi^*} \bar{\Gamma}^k_{ij}=\Gamma_{ij}^k$.

%\label{Christoffel_Isometry}\end{equation}
%\begin{Prop} Let $M$ and $N$ be connected \marginpar{complete?!}Riemannian supermanifolds and $\Phi:M\to N$ an affine mapping. If $d_p\Phi:T_pM\to T_{\phi(p)}N$ is an isometry for some point $p$, then $\Phi$ is an isometry. 
%\end{Prop}
%\begin{proof} Pick some point $q\in M$ and some supergeodesic $\gamma$ joining $p$ and $q$. 
%\end{proof}
The following lemma is clear.
\begin{Lemma} \label{Isometry_Geodesic}Let $\gamma:\R^{1|1}\to M$ be a supergeodesic and $\Phi:M\to N$ an isometry. Then $\Phi\circ \gamma$ is a supergeodesic.
\end{Lemma}

\begin{Lemma}\label{LemmaGeodesics} Let $M$ be a supermanifold, $(x_i,\xi_\alpha)$ local coordinates on $U\subset M$ and $p$ some point of $U$. Let furthermore $f$ be a function on $U$ of degree one with respect to these coordinates. If for all supergeodesics $\gamma$ starting at $p$ we have $\gamma^*f=0$, then $f=0$ on some neighbourhood of $p$.
\end{Lemma}
\begin{proof} Write $f$ in these coordinates: $f=\sum_\alpha f_\alpha\xi_\alpha$ for some $\shC^\infty$-functions $f_\alpha$. Assume that $f$ is not zero in any neighbourhood of $p$; then we can find some $q$ near $p$ such that there exists an (ordinary) geodesic $\widetilde{\gamma}$ with $\widetilde{\gamma}(0)=p$, $\widetilde{\gamma}(1)=q$ and such that not all $f_\alpha(q)=0$. Let $\gamma$ be the supergeodesic with underlying curve $\widetilde{\gamma}$ that satisfies the following condition: If $\gamma^*\xi_\alpha=h_\alpha\cdot \xi$ for some smooth functions $h_\alpha$, then $h_\alpha(1)=f_\alpha(q)$ (recall that the odd differential equations (\ref{RightOddGeodesicEquation}) are of first order). Then we have $$(\gamma^*f)[1]=\sum_\alpha f_\alpha(\widetilde{\gamma}(1)) (\gamma^*\xi_\alpha)[1]=\sum_\alpha f_\alpha^2(q)\xi\neq 0,$$ which contradicts the assumption. Here, for a superfunction \mbox{$g=g_0+ g_1\xi$} on $\R^{1|1}$ written in coordinates, $g[1]$ shall denote the element of $\Lambda_\R[\xi]$ one gets by inserting $1$ into the coefficient functions, i.e.~$g[1]=g_0(1)+ g_1(1)\xi$.
\end{proof}

%\begin{Lemma}\label{LemmaGeodesics} Let $M$ be a supermanifold and $f$ a local superfunction on $M$ of degree one with respect to coordinates $(x_i,\xi_\alpha)$ around $p\in M$. If for all supergeodesics $\gamma$ starting at $p$, $\gamma^*f=0$, then $f=0$ on some neighbourhood of $p$.
%\end{Lemma}
%\begin{proof} Write $f$ in these coordinates: $f=\sum_\alpha f_\alpha\xi_\alpha$ for some $\shC^\infty$-functions $f_\alpha$. Assume that $f$ is not zero in any neighbourhood of $p$; then we can find some $q$ near $p$ such that there exists an (ordinary) geodesic $\widetilde{\gamma}$ with $\widetilde{\gamma}(0)=p$, $\widetilde{\gamma}(1)=q$ and such that not all $f_\alpha(q)=0$. We can lift $\widetilde{\gamma}$ horizontally to a curve $\bar{\gamma}$ such that $\bar{\gamma}(1)=\sum_\alpha f_\alpha(q)\xi_\alpha(q)$ -- here we regard the $\xi_\alpha$ as sections of $E$ so that the expression $\xi_\alpha(q)$ makes sense. Then for the associated supergeodesic $\gamma$, $$(\gamma^*f)[1]=\sum_\alpha f_\alpha(\widetilde{\gamma}(1)) (\gamma^*\xi_\alpha)[1]=\sum_\alpha f_\alpha^2(q)\xi\neq 0;$$ thus, in particular, $\gamma^*f\neq 0$. Here, for a superfunction $g=g_0+ g_1\xi$ on $\R^{1|1}$ written in coordinates, $g[1]$ shall denote the element of $\Lambda_\R[\xi]$ one gets by inserting $1$ into the coefficient functions, i.e.~$g[1]=g_0(1)+ g_1(1)\xi$.
%\end{proof}
\begin{Rem} For superfunctions of degree zero, i.e.~$\shC^\infty$-functions, the analogous statement follows directly from the fact that the exponential map is a local diffeomorphism. For superfunctions of higher degree, the lemma is false.
\end{Rem}
Before we are able to prove that isometries are determined by the data at one point (just as geodesics, cf. Proposition \ref{GeodesicDetermined}), we need a technical lemma:
\begin{Lemma} \label{TechnicalLemmaCross}Let $M$ be a supermanifold with odd dimension $q$ and $(x_i,\xi_\alpha)$ coordinates on $U\subset M$. If $f_1,\ldots,f_q$ are superfunctions on $U$ such that $$\partial_\alpha f_\beta=\partial_\beta f_\alpha$$ for all $\alpha,\beta$, then all the $f_\alpha$ are the sum of functions of degree at most one (with respect to the chosen coordinates).
\end{Lemma}
\begin{proof} Assume the contrary, i.e.~at least one of the functions contains some summand of degree higher than one, say $f_1$. Then there are $1\le \alpha,\beta\le q$ such that $$f_1=\xi_\alpha\xi_\beta g_1+\xi_\alpha g_2+\xi_\beta g_3 + g_4$$ for superfunctions $g_1,\ldots,g_4$ not containing $\xi_\alpha$ and $\xi_\beta$ with at least $g_1\neq 0$. Then, \begin{equation}\label{g_1KeinXi_1}\partial_1 f_\alpha = \partial_\alpha f_1=\xi_\beta g_1 +g_2,\end{equation} and thus $$f_\alpha=\xi_1(\xi_\beta g_1+g_2)+g_5$$ for some function $g_5$ not containing $\xi_1$. Note also that we see from (\ref{g_1KeinXi_1}) that $g_1$ does not contain $\xi_1$ since both $g_1$ and $g_2$ do not contain $\xi_\beta$. Analogously, $$\partial_1 f_\beta=\partial_\beta f_1=-\xi_\alpha g_1+g_3$$ implies $$f_\beta=\xi_1(-\xi_\alpha g_1+g_3)+g_6$$ for some $g_6$ not containing $\xi_1$. But then on the one hand
$$\partial_\beta f_\alpha=-\xi_1 g_1+\partial_\beta g_5$$ and on the other hand $$\partial_\beta f_\alpha=\partial_\alpha f_\beta=\xi_1 g_1+\partial_\alpha g_6$$ which is a contradiction since $g_1\neq 0$ and $g_1,g_5$ and $g_6$ do not contain $\xi_1$.
\end{proof}

\begin{Prop} \label{IsometryDetermined}An isometry of a connected Riemannian supermanifold $M$ is determined by its value and its derivative at one point.
\end{Prop}
\begin{proof} Let $\Phi=(\phi,\phi^*):M\to M$ be an isometry and $p\in M_{\red}$ such that $\phi(p)=p$ and $d_p\Phi=\id_{T_pM}$. We have to show $\Phi=\id_M$. The corresponding result in the standard theory says that $\phi=\id_{M_{\red}}$.

Fix coordinates on some open set $U$ containing $p$ and write %$$\gamma^*\xi_\alpha=h_\alpha \xi$$ 
$$\phi^*\xi_\alpha=\sum_{\beta} f^\alpha_\beta \xi_\beta +O(\xi^3)$$ for some $\shC^\infty$-functions $f^\alpha_\beta$. For every supergeodesic $\gamma$ starting at $p$, we have $\Phi\circ \gamma=\gamma$ because of Lemma \ref{Isometry_Geodesic} and Proposition \ref{GeodesicDetermined}. In coordinates, this means $$\gamma^*\left(\xi_\alpha-\sum_\beta f^\alpha_\beta \xi_\beta\right)=0$$ for all $\alpha$. Lemma \ref{LemmaGeodesics} now says that $\sum_\beta f^\alpha_\beta \xi_\beta=\xi_\alpha$, probably after restricting $U$. But since the $\xi_\alpha$ are linearly independent as elements of the $\shC^\infty(U)$-module $\shO_M(U)$, $f^\alpha_\beta=\delta_{\alpha\beta}$. Summing up, we have shown $$\phi^*x_i=x_i+O(\xi^2),\quad \phi^*\xi_\alpha=\xi_\alpha+O(\xi^3).$$
%Consequently, $d\Phi(\partial_\alpha)=\partial_\alpha+N_\alpha$, where $N_\alpha$ is an odd vector field with only nilpotent coefficients, i.e.~the value of $N_\alpha$ at any point is $0$. Since $\phi=\Id_{M_{\red}}$ implies $\phi^*x_i=x_i+O(\xi^2)$, we see that $d\Phi(\partial_i)=\partial_i+N_i$, where $N_i$ is an even nilpotent vector field. 
%Hence, for any vector field $X$, $d\Phi(X)=X+N$ for some nilpotent vector field $N$, since the values of $X$ and $d\Phi(X)$ coincide at each point. Now we may express $d\Phi$ in terms of a local orthonormal frame field $(X_i,Y_\alpha)$ (cf. Proposition \ref{PropONFrame}) as follows:
%\begin{align}d\Phi(X_i)=X_i+\sum_j f^i_j X_j+\sum_\beta f^i_\beta Y_\beta,\\
%d\Phi(Y_\alpha)=Y_\alpha+\sum_j f^\alpha_j X_j+\sum_\beta f^\alpha_\beta Y_\beta,\end{align} where all appearing coefficient functions are nilpotent. \marginpar{Beziehungen zwischen $f$'s???}
We continue by showing that the terms of higher degree vanish, starting with the smallest: Writing %$$\phi^*x_i=x_i+\sum_{\alpha_1<\alpha_2} f^i_{\alpha_1\alpha_2}\xi_{\alpha_1}\xi_{\alpha_2}+O(\xi^4),$$
$$\phi^*x_i=x_i+f^i_2+O(\xi^4),$$ where $f^i_2$ is homogeneous of degree $2$ with respect to the $\Z$-grading, we calculate
%\begin{align*}d\Phi(\partial_i)&=\partial_i+\sum_j\left(\sum_{\alpha_1<\alpha_2}(\partial_i f^j_{\alpha_1\alpha_2})\xi_{\alpha_1}\xi_{\alpha_2}\right)\partial_j+O(\xi^3),\\
%d\Phi(\partial_\alpha)&=\partial_\alpha+\sum_j\left(\sum_{\alpha<\alpha_2} f^j_{\alpha\alpha_2}\xi_{\alpha_2}-\sum_{\alpha_1<\alpha}f^j_{\alpha_1\alpha}\xi_{\alpha_1}\right)\partial_j+O(\xi^2),\end{align*}
\begin{align*}d\Phi(\partial_i)&=\partial_i+\sum_j(\partial_i f^j_2)\partial_j+O(\xi^3),\quad
d\Phi(\partial_\alpha)=\partial_\alpha+\sum_j(\partial_\alpha f^j_2)\partial_j+O(\xi^2),\end{align*} 
where the $O$-notation is to be understood for the coefficient functions, when the vector fields are expressed in the basis $\partial_i,\partial_\alpha$. Since $\Phi$ is an isometry, \begin{align*}\left<\partial_i,\partial_\alpha\right>&=\phi^*\left<d\Phi(\partial_i),d\Phi(\partial_\alpha)\right>
=\left<\partial_i,\partial_\alpha\right>+\sum_j(\partial_\alpha f^j_2)\left<\partial_i,\partial_j\right>+O(\xi^2).\end{align*} Thus, $$\sum_j (\partial_\alpha f_2^j)\left<\partial_i,\partial_j\right>^{\widetilde{}}=0.$$ The matrix $(\left<\partial_i,\partial_j\right>^{\widetilde{}}\,)_{i,j}$ is invertible, so we get $\partial_\alpha f^i_2=0$ for all $\alpha$. Since $f^i_2$ is homogeneous of degree $2$, this is only possible if $f^i_2=0$. Thus, $$\phi^*x_i=x_i+O(\xi^4).$$
Now we deal with the odd coordinates:
%$$\phi^*\xi_\alpha=\xi_\alpha+\sum_{\alpha_1<\alpha_2<\alpha_3}f^\alpha_{\alpha_1\alpha_2\alpha_3}\xi_{\alpha_1}\xi_{\alpha_2}\xi_{\alpha_3}+O(\xi^5)$$
From $$\phi^*\xi_\alpha=\xi_\alpha+f^\alpha_3+O(\xi^5)$$ for some homogeneous functions $f^\alpha_3$ of degree $3$ we get
%\begin{align*} d\Phi(\partial_i)&=\partial_i+\sum_\beta\left(\sum_{\alpha_1<\alpha_2<\alpha_3}(\partial_i f^\beta_{\alpha_1\alpha_2\alpha_3})\xi_{\alpha_1}\xi_{\alpha_2}\xi_{\alpha_3}\right)\partial_\beta+O(\xi^4),\\
%d\Phi(\partial_\alpha)&=\partial_\alpha+\sum_\beta\left(\sum_{\alpha<\alpha_2<\alpha_3} f^\beta_{\alpha\alpha_2\alpha_3}\xi_{\alpha_2}\xi_{\alpha_3}-\sum_{\alpha_1<\alpha<\alpha_3} f^\beta_{\alpha_1\alpha\alpha_3}\xi_{\alpha_1}\xi_{\alpha_3}\right.\\ &\quad\quad\quad\quad\quad\quad\quad\left.+\sum_{\alpha_1<\alpha_2<\alpha}f^\beta_{\alpha_1\alpha_2\alpha}\xi_{\alpha_1}\xi_{\alpha_2}\right)\partial_\beta+O(\xi^3)\\
%&=: \partial_\alpha+\sum_\beta g^\beta_\alpha\partial_\beta+O(\xi^3),\end{align*}
%where we introduced the functions $g^\beta_\alpha$ (which are of degree two) for brevity. 
\begin{align*} %d\Phi(\partial_i)&=\partial_i+\sum_\beta (\partial_i f^\beta_3)\partial_\beta+O(\xi^4),\\
d\Phi(\partial_\alpha)&=\partial_\alpha+\sum_\beta (\partial_\alpha f^\beta_3)\partial_\beta+O(\xi^3).
\end{align*}
Then
$$\left<\partial_\alpha,\partial_\delta\right>=\left<\partial_\alpha,\partial_\delta\right>+\sum_\beta \left((\partial_\delta f^\beta_3)\left<\partial_\alpha,\partial_\beta\right>-(\partial_\alpha f^\beta_3) \left<\partial_\delta,\partial_\beta\right>\right)+O(\xi^3),$$ which implies \begin{align*}\partial_\delta\left(\sum_\beta f^\beta_3 \left<\partial_\alpha,\partial_\beta\right>^{\widetilde{}}\right)&=\sum_\beta (\partial_\delta f^\beta_3) \left<\partial_\alpha,\partial_\beta\right>^{\widetilde{}}\\
&=\sum_\beta(\partial_\alpha f^\beta_3) \left<\partial_\delta,\partial_\beta\right>^{\widetilde{}}=\partial_\alpha\left(\sum_\beta f^\beta_3 \left<\partial_\delta,\partial_\beta\right>^{\widetilde{}}\right).\end{align*} But now Lemma \ref{TechnicalLemmaCross} shows that $\sum_\beta f^\beta_3 \left<\partial_\alpha,\partial_\beta\right>^{\widetilde{}}=0$ for all $\alpha$. Since the matrix $(\left<\partial_\alpha,\partial_\beta\right>^{\widetilde{}}\,)_{\alpha,\beta}$ is invertible, all the $f^\alpha_3=0$; we have shown $$\phi^*\xi_\alpha=\xi_\alpha+O(\xi^5).$$ It is clear that 
we can proceed inductively on the degree the same way, alternatingly dealing with the even and odd coordinates -- we only used that the $f^i_2$ and $f^\alpha_3$ are homogeneous of degree greater than one. Then we see that $\Phi$ is the identity on $U$; an easy argument using the connectedness of $M$ now proves that $\Phi$ is the identity on the whole of $M$.\end{proof}

\subsection{Graded Killing Fields}

Before attacking the problem of introducing more structure to the isometry group of a Riemannian supermanifold, we deal with the corresponding infinitesimal objects: the Killing vector fields.

Let $M$ be a Riemannian supermanifold. A {\it graded Killing vector field on} $M$ is a vector field $X$ such that
\begin{equation}
X\left<Y,Z\right>=\left<[X,Y],Z\right>+(-1)^{|X||Y|}\left<Y,[X,Z]\right>
\end{equation} for all $Y,Z$. Thus, using the properties of the Levi-Civita connection $\nabla$, $X$ is a graded Killing vector field if and only if
%%%%%%%%% IM ARTIKEL WEG
%\begin{align*}
%\left<\nabla_XY,Z\right>+(-1)^{|X||Y|}\left<Y,\nabla_XZ\right> =& \left<\nabla_XY-(-1)^{|X||Y|}\nabla_YX,Z\right>\\ &+(-1)^{|X||Y|}\left<Y,\nabla_XZ-(-1)^{|X||Z|}\nabla_ZX\right>,
%\end{align*} i.e.~
%%%%%%%%% BIS HIER
\begin{equation} \left<\nabla_YX,Z\right>+(-1)^{|X||Y|+|X||Z|+|Y||Z|}\left<\nabla_ZX,Y\right>=0. \label{Killing_Equation}
\end{equation}

%Defining the {\it Lie derivative} ${\mathcal L}_X$ with respect to an even vector field $X$ as in \cite{DelMor}, p. 79, and using the elementary properties proven there, we see that if $X$ is an even vector field, it is Killing if and only if
%${\mathcal L}_X\left<\cdot,\cdot\right>=0.$

The super vector space of all graded Killing vector fields becomes a Lie superalgebra, if we define the bracket to be induced from the Lie superalgebra of all vector fields. Following usual conventions, we will denote by {\it the Lie superalgebra of graded Killing vector fields} the opposite to this Lie superalgebra, i.e.~the Lie superalgebra 
with the same underlying vector space and the new bracket the negative of the old bracket. 

We denote the second covariant derivative by $\nabla^2_{Y,Z}X:=\nabla_Y\nabla_Z X-\nabla_{\nabla_Y Z}X$. Then we have \begin{align}R(Y,Z)X
%%%%%%%%%%%%% IM ARTIKEL WEG
%&=\nabla_Y\nabla_ZX -(-1)^{|Y||Z|}\nabla_Z\nabla_YX-\nabla_{[Y,Z]}X\nonumber\\ &=\nabla_Y\nabla_Z X-\nabla_{\nabla_YZ}X-(-1)^{|Y||Z|}(\nabla_Z\nabla_YX-\nabla_{\nabla_ZY}X)\nonumber\\
%%%%%%%%%%%%% BIS HIER
&=\nabla^2_{Y,Z}X-(-1)^{|Y||Z|}\nabla^2_{Z,Y}X.\label{Curv_2ndCovDer}\end{align} 
\begin{Lemma}\label{Killing_2ndCovDer_Skew} Let $X$ be a Killing vector field. Then $$\left<\nabla^2_{Y,Z}X,W\right>+(-1)^{|Z||X|+|Z||W|+|X||W|}\left<\nabla^2_{Y,W}X,Z\right>=0.$$
\end{Lemma}
\begin{proof} This is a direct calculation using (\ref{Killing_Equation}) several times.
%%%%%%%%%%%%%%%%% IM ARTIKEL WEG
% two times:
%\begin{align*} \left<\nabla^2_{Y,Z}X,W\right>&=\left<\nabla_Y\nabla_Z X-\nabla_{\nabla_Y Z}X,W\right> \\
%&=\left<\nabla_Y\nabla_Z X,W\right>+(-1)^{(|Y|+|Z|)(|X|+|W|)+|X||W|}\left<\nabla_W X,\nabla_Y Z\right>\\
%&=Y\left<\nabla_Z X,W\right>-(-1)^{|Y|(|Z|+|X|)}\left<\nabla_Z X,\nabla_Y W\right>\\
%&\quad +(-1)^{|Z||X|+|Z||W|+|X||W|}(Y\left<\nabla_W X,Z\right>-\left<\nabla_Y\nabla_W X,Z\right>)\\
%&=-(-1)^{|Z||X|+|Z||W|+|X||W|}\left<\nabla^2_{Y,W}X,Z\right>.
%\end{align*}
%%%%%%%%%%%%%%%%% BIS HIER
\end{proof}

\begin{Prop}\label{KillingDGL}Let $X$ be a Killing vector field. Then $$\nabla^2_{Y,Z}X=-(-1)^{|X|(|Y|+|Z|)}R(X,Y)Z.$$
\end{Prop}
\begin{proof} The proof is analogously to \cite{Petersen}, p.216, with a lot of additional signs; it uses Lemma \ref{Killing_2ndCovDer_Skew}, formula (\ref{Curv_2ndCovDer}) and the symmetries of the curvature tensor collected in Proposition \ref{CurvProperties}. Since no new idea enters, we again omit the calculation.
\end{proof}

\begin{Prop} Let $X$ be a Killing vector field on a Riemannian supermanifold. If there exists a point $p$ such that $X(p)=0$ and $(\nabla X)(p)=0$, then $X=0$. \label{KillingDeterminedByPoint}
\end{Prop}
\begin{Rem} In standard Riemannian geometry, this is proven either via the flow of $X$ or the fact that $X$, restricted to any geodesic, is a Jacobi field. Since we did not introduce Jacobi fields and flows of odd vector fields, we need a different proof. 
\end{Rem}
\begin{proof}
Consider first the case of $X$ even. Pick coordinates $(x_i,\xi_\alpha)$ around $p$ and write $X=\sum_i (\sum_j f^i_j)\partial_i$ in coordinates; we assume that the degree of $f^i_j$ is $j$. Proposition \ref{KillingDGL}, written in coordinates, is a system of second order differential equations; it yields the vanishing of the coefficient functions of degree $0$ and $1$, i.e.~$$X=\sum_i f^i_2\partial_i+O(\xi^3).$$
(Recall that the $O$-notation is to be understood for the coefficient functions.) Then on the one hand $X\left<\partial_j,\partial_\beta\right>=O(\xi^3)$ and on the other hand \begin{align*}X\left<\partial_j,\partial_\beta\right>&=\left<[X,\partial_j],\partial_\beta\right>+\left<\partial_j,[X,\partial_\beta]\right>
=-\sum_i (\partial_\beta f^i_2)\left<\partial_j,\partial_i\right>^{\widetilde{}}+O(\xi^3).\end{align*} Since the matrix $(\left<\partial_i,\partial_j\right>^{\widetilde{}})$ is invertible, we conclude $\partial_\beta f^i_2=0$ for all $i$ and $\beta$. Thus, $f^i_2=0$, since the function was assumed to be homogeneous of positive degree.

Attacking the remaining term of smallest degree we write $$X=\sum_\alpha f^\alpha_3\partial_\alpha+O(\xi^4).$$ Then on the one hand $X\left<\partial_\beta,\partial_\delta\right>=O(\xi^4)$ and on the other hand \begin{align*}X\left<\partial_\beta,\partial_\delta\right>&=\left<[X,\partial_\beta],\partial_\delta\right>+\left<\partial_\beta,[X,\partial_\delta]\right>\\ &=-\sum_\alpha(\partial_\beta f^\alpha_3)\left<\partial_\alpha,\partial_\delta\right>^{\widetilde{}}-\sum_\alpha (\partial_\delta f^\alpha_3)\left<\partial_\beta,\partial_\alpha\right>^{\widetilde{}}+O(\xi^4).\end{align*} Lemma \ref{TechnicalLemmaCross} yields the vanishing of the $f^\alpha_3$ via the same argument as in the last part of the proof of Proposition \ref{IsometryDetermined}; then induction finishes the case of $X$ even.

If $X$ is odd, we can argue very similarly, with only some signs changing. First of all, the same argument shows the vanishing of the coefficient functions of degree $0$ and $1$. Applying $X=\sum_\alpha f^\alpha_{2}\partial_\alpha+O(\xi^{3})$ to the superfunctions $\left<\partial_\beta,\partial_\delta\right>$ yields the vanishing of the $f^\alpha_{2}$; thereafter applying $X=\sum_i f^i_{3}\partial_i+O(\xi^{4})$ to $\left<\partial_j,\partial_\beta\right>$ shows that $f^i_{3}=0$. As usual we continue by induction.
\end{proof}

\subsection{The Isometry Group}

In this section we want to define the isometry group $I(M)$ of a Riemannian supermanifold $M$; naturally, it shall become not a Lie group but a Lie supergroup. Staying consistent with the ungraded case, the Lie superalgebra of $I(M)$ has to be the Lie superalgebra of graded Killing vector fields. 

The set of all isometries (cf.  \ref{SectionIsometries}) of a Riemannian supermanifold $M$ clearly is a group. Our task now is to show that there is a natural Lie structure on it. Then we will equip this Lie group with the Lie superalgebra of graded Killing vector fields to get a Harish-Chandra pair  -- this will by definition give us the isometry supergroup of $M$.
 
In spirit of this, we denote by $I(M)_{\red}$ the group of isometries of $M$, although we have not yet defined $I(M)$. This is not to be mixed up with $I(M_{\red})$, which is simply the usual isometry group of the pseudo-Riemannian manifold $M_{\red}$. Recall that the isometry group of a conventional pseudo-Riemannian manifold, endowed with the compact-open topology, is a Lie group, see e.g.~Ballmann \cite{Ballmann}: it is a closed subgroup of the group of affine diffeomorphisms with respect to the Levi-Civita connection of the metric.

It would be nice to say that we equip $I(M)_{\red}$ with the compact-open topology; but since topological notions like compactness and openness make sense only for the underlying manifold, this is not possible. In the following we therefore want to realize $I(M)_{\red}$ as a closed subgroup of an automorphism group of the parallelization of some ordinary differentiable manifold and then use the general result of Ballmann, namely that such groups are Lie groups. 

Let us briefly recall the definitions of \cite{Ballmann}: If $N$ is a conventional $n$-dimensional manifold together with a parallelization $\Phi: N\times \R^n\to TN$, any $z\in \R^n$ induces a so-called {\it constant vector field} $Z(p)= \Phi(p,z)$. The {\it automorphism group} $\Aut(\Phi)$ of $\Phi$ then is the group of all diffeomorphisms \mbox{$f:N\to N$} such that $df\circ Z=Z\circ f$ for all constant vector fields $Z$, equipped with the compact-open topology.

The following is an extension of Example $3.1$ of \cite{Ballmann}: Let $M$ be an $n|2m$-dimensional Riemannian supermanifold with Levi-Civita connection $\nabla$ and consider the vector bundle $TM=(TM)_0\oplus (TM)_1\to M_{\red}$ with \mbox{$n+2m$-dimensional} fibres $T_pM$ (we may regard it as a bundle of $n|2m$-dimensio\-nal super vector spaces). Recall that the connection on $M$ induces a connection on $TM\to M_{\red}$. Let $\pi:\GL(TM)\to M_{\red}$ be the $\GL(n)\times \GL(2m)$-principal fibre bundle of graded frames; the fibre over $p\in M_{\red}$ is given by $$\GL(TM)_p=\GL(\R^n,(T_pM)_0)\times\GL(\R^{2m},(T_pM)_1).$$ We want to find a parallelization of the conventional manifold $\GL(TM)$ such that affine diffeomorphisms of the supermanifold $M$ induce diffeomorphisms of $\GL(TM)$ respecting the parallelization.

The vertical distribution ${\mathcal V}$ of $\GL(TM)$ given by the kernel of $d\pi$ is trivialized via the mapping $$\Phi_v:\GL(TM)\times (\gl(n)\times \gl(2m))\to {\mathcal V}$$ defined by $$\Phi_v((\phi_1,\phi_2),(x_1,x_2))=\partial_t\left.(\phi_1\exp(tx_1),\phi_2\exp(tx_2))\right|_{t=0}.$$ It remains to find a corresponding (trivializable) horizontal distribution ${\mathcal H}$ (which will be of rank $n$). For any $\phi=(\phi_1,\phi_2)\in \GL(TM)$ and $z\in \R^n$, let $v=\phi_1(z)\in T_pM_{\red}$, where $p=\pi(\phi)\in M_{\red}$. Choose a curve $c$ in $M_{\red}$ through $p$ with $c'(0)=v$ and let $\psi$ be the unique parallel frame along $c$ such that $\psi(0)=\phi$. Then $$\Phi_h(\phi,z):=\psi'(0)$$ gives a trivialization $$\Phi_h:\GL(TM)\times \R^n\to {\mathcal H}$$ of a horizontal distribution ${\mathcal H}$. The pair $\Phi=(\Phi_v,\Phi_h)$ therefore is a parallelization of $\GL(TM)$. If $f$ is a diffeomorphism of the supermanifold $M$, it induces a diffeomorphism $f_*$ of $\GL(TM)$ sending a frame $\phi$ to $df\circ \phi$; if $f$ is affine, this diffeomorphism is an element of $\Aut(\Phi)$.

We thus get a mapping $I(M)_{\red}\to \Aut(\Phi)$, which is injective by Proposition \ref{IsometryDetermined}. 
\begin{Prop} The image of the natural inclusion $I(M)_{\red}\to \Aut(\Phi)$ is closed in $\Aut(\Phi)$.
\end{Prop}
\begin{proof} Let $g_n:M\to M$ be a sequence of isometries of the supermanifold $M$ such that the induced mappings ${g_{n}}_*$ converge against some $h$ with respect to the compact-open topology.

The space $\shO_M(M)$ is a Fr\'echet space via the family of semi-norms $\vert\cdot\vert_{K,\partial}$ defined as follows: for any compact subset $K\subset M_{\red}$ and any differential operator $\partial$, we set $\vert f\vert_{K,\partial}:=\sup_{p\in K}\vert (\partial f)(p)\vert$, see Kostant \cite{Kostant}, p.199. 
By writing the condition of being an isometry in coordinates, an isometry can be viewed as a solution of some system of ordinary differential equations. Since the ${g_n}_*$ converge, also these solutions converge in the sense that for any $f\in \shO_M(M)$, the sequence $g_n^*f$ converges against some $\Gamma(f)\in\shO_M(M)$. Then, the mapping $\Gamma:\shO_M(M)\to \shO_M(M)$ defines an isometry $g$ of $M$ such that the induced mapping $g_*:\GL(TM)\to \GL(TM)$ equals $h$.
\end{proof}
Thus, by giving $I(M)_{\red}$ the induced topology of $\Aut(\Phi)$, the results of \cite{Ballmann} introduce the structure of Lie group on $I(M)_{\red}$.

\begin{Rem} Another possible way of turning $I(M)_{\red}$ into a Lie group is the following: Give it the coarsest topology such that for all $f\in \shO_M(M)$, the mapping $I(M)_{\red}\to \shO_M(M);\, g\mapsto g^*f$ is continuous -- here, $\shO_M(M)$ carries the structure of a Fr\'echet space mentioned in the proof above. Then, $I(M)_{\red}$ is a locally compact topological group without small subgroups. It is known that such groups are Lie groups, see \cite{MontgomeryZippin}, p.107 and the references given there.
\end{Rem}

To complete the construction of the isometry supergroup we have to equip this Lie group with a suitable Lie superalgebra.

Let us first calculate the Lie algebra of $I(M)_{\red}$. If we take a left invariant vector field $X$ on $I(M)_{\red}$, it defines a unique one-parameter group $(g_t)$ in $I(M)_{\red}$. Then the associated vector field $f\mapsto \left.\partial_t\right|_{t=0}(g_t^*f)$ is an even Killing field; this correspondence is an isomorphism from the Lie algebra of $I(M)_{\red}$ to the Lie algebra of even Killing fields on $M$. %recall the sign convention in the definition of the Lie superalgebra of graded Killing vector fields. \marginpar{Lieklammer nachrechnen}

The {\it isometry group} $I(M)$ of $M$ is now the Lie supergroup associated to the Harish-Chandra pair $$(I(M)_{\red},\frak g),$$ where $\frak g$ is the Lie superalgebra of all graded Killing vector fields on $M$. We have already verified that the Lie algebra of $I(M)_{\red}$ is equal to $\mathfrak g_0$, so what is left is to specify the action of $I(M)_{\red}$ on $\mathfrak g$. %One is tempted to simply say that it is the usual adjoint action but remember that we did not define flows for odd vector fields -- there is no correspondence between odd Killing vector fields on $M$ and vector fields on $I(M)_{\red}$. 
For a graded Killing vector field $X$ on $M$ and an isometry $g\in I(M)_{\red}$, we define $$\Ad_g X:=dg(X)=(g^{-1})^*\circ X\circ g^*.$$ On $\mathfrak g_0$, this coincides with the action of $I(M)_{\red}$ on its Lie algebra. Note that if $M$ is a usual manifold, $I(M)$ is the usual isometry group. 

There is a natural action of $I(M)$ on $M$ which we give in terms of the above Harish-Chandra pair, cf. \ref{Section_Action}: The Lie group $I(M)_{\red}$ clearly acts on $M$; in fact, it is defined as a group of diffeomorphisms of $M$. We have to define a compatible morphism from the Lie superalgebra of graded Killing vector fields on $M$ to the opposite of the Lie superalgebra of vector fields on $M$. Here, we simply take the inclusion.

\subsection{Invariant Metrics on Lie Supergroups}

%Different from the standard theory, it will turn out that the super Killing form (which will be introduced in \ref{SectionKilling}) will only in the least cases have a definite sign. Recall that therefore, our notion of Riemannian metric on supermanifolds did not include any positivity condition.

Given an ordinary Lie group $G$, a pseudo-Riemannian metric $\left<\cdot,\cdot\right>$ on $G$ is called {\it left-invariant} if the left translations are isometries. Equivalent to this is the condition that for all left-invariant vector fields $X$ and $Y$, $\left<X,Y\right>$ is a constant function on $G$ (i.e.~a real multiple of the unit in the algebra of global functions on $G$). We take this to be the definition of left-invariance in the case of Lie supergroups. Analogously, we define right-invariance and then bi-invariance.

The left-invariant graded Riemannian metrics on a Lie supergroup $G$ are in one-to-one correspondence with the scalar superproducts on ${\mathfrak g}$: for $X,Y\in {\mathfrak g}$ regard the real number $\left< X,Y\right>_e$ as a constant function on $G$ and extend the metric linearly with respect to superfunctions.

A scalar superproduct $\left<\cdot,\cdot\right>_e$ on $\mathfrak g$ is called {\it $\Ad_{G_{\red}}$-invariant} if $$\left<\Ad_gX,\Ad_gY\right>_e=\left<X,Y\right>_e$$ for all $X,Y\in \mathfrak g$ and all $g \in G_{\red}$. It is {\it $\ad_{\mathfrak{g}}$-invariant} if $$\left<[X,Y],Z\right>_e+(-1)^{|X||Y|}\left<Y,[X,Z]\right>_e=0$$ for all $X,Y,Z\in \mathfrak g$. Of course, if $G$ is connected, $\ad_{\mathfrak g}$-invariance implies $\Ad_{G_{\red}}$-invariance. We call the scalar superproduct {\it $\Ad_G$-invariant} or simply {\it $\Ad$-invariant} if it is $\Ad_{G_{\red}}$- and $\ad_{\mathfrak{g}}$-invariant.

%%%%% IM ARTIKEL DAZU

As an example, $G=\R^{1|2}$ with the Lie supergroup structure given in symbolic notation by $$(x,\xi_1,\xi_2)+(t,\theta_1,\theta_2):=(x+t+\xi_1\xi_2+\theta_1\theta_2,\xi_1+\theta_1,\xi_2+\theta_2)$$ does not admit any bi-invariant metric; its Lie algebra admits $\Ad_{G_{\red}}$-invariant but no $\ad_{\mathfrak g}$-invariant scalar superproducts.

%%%%%%%%%%% IM ARTIKEL WEG
%If a left-invariant graded Riemannian metric $\left<\cdot,\cdot\right>$ on a connected Lie supergroup $G$ is given, $\ad_{\mathfrak{g}}$-invariance implies $\Ad_{G_{\red}}$-invariance. The converse is not true -- which is reasonable since $\Ad_{G_{\red}}$ is only a representation of $G_{\red}$ -- as the following example shows:  Consider on $G=\R^{1|2}$ the Lie supergroup structure given in symbolic notation by $$(x,\xi_1,\xi_2)+(t,\theta_1,\theta_2):=(x+t+\xi_1\xi_2+\theta_1\theta_2,\xi_1+\theta_1,\xi_2+\theta_2).$$
%Then an easy calculation shows that the left-invariant vector fields on $G$ are spanned by $$X_1=\partial_x,\quad X_2=-\xi_1\partial_x+\partial_{\xi_1},\quad X_3=-\xi_2\partial_x+\partial_{\xi_2},$$ whereas the right-invariant vector fields are spanned by
%$$Y_1=\partial_x,\quad Y_2=\xi_1\partial_x+\partial_{\xi_1},\quad Y_3=\xi_2\partial_x+\partial_{\xi_2}.$$
%The left-invariant metric defined by $$\left<X_1,X_1\right>=1,\quad \left<X_2,X_3\right>=1$$ now induces an $\Ad_{G_{\red}}$-invariant scalar superproduct at $e$, which is not $\ad_{\mathfrak{g}}$-invariant (regard e.g.~$\left<[X_2,X_2],X_1\right>-\left<X_2,[X_2,X_1]\right>\neq 0$). The metric is not bi-invariant since e.g.~$\left<Y_1,Y_2\right>=2\xi_1\neq 0$; there even does not exist any bi-invariant metric on $G$.
%%%%%%%%%%% BIS HIER

The equivalence of the bi-invariance of a metric and the $\ad$-invariance of the induced supersymmetric bilinear form at $e$, which is very easily proven in the standard case, requires some preparation before we can establish it in Theorem \ref{BiAd}. 
\begin{Lemma} Let $G$ be a Lie supergroup. Then for any left-invariant vector field $X$ and any right-invariant vector field $Y$, we have $[X,Y]=0$.\label{BracketLeftRight}
\end{Lemma}
\begin{proof} Let $\tau,\sigma\in T_eG$ and $X_\tau,X_\sigma$ be the corresponding left-invariant vector fields. The differential of the inverse map $i$ interchanges left- and right-invariant vector fields, so it suffices to show $[X_\tau,di(X_\sigma)]=0$; this is an easy calculation using (\ref{VectorToLInvVF}).
%%%%%%%%%%% IM ARTIKEL WEG
%\begin{align*} X_\tau\circ di(X_\sigma) (f)&= (I\otimes \tau)\circ m^*\circ i^*\circ (I\otimes \sigma)\circ m^*\circ i^*(f)\\
%&=(I\otimes \tau)\circ m^*\circ i^*\circ (I\otimes \sigma)(\{(g,h)\mapsto f(h^{-1}g^{-1})\})\\
%&=(I\otimes \tau)\circ m^*(\{g\mapsto \sigma(\{h\mapsto f(h^{-1}g)\})\})\\
%&=\{g\mapsto \tau(\{k\mapsto \sigma(\{h\mapsto f(h^{-1}gk)\})\})\},\\ \\
%di(X_\sigma)\circ X_\tau(f)&=i^*\circ (I\otimes \sigma)\circ m^*\circ i^*\circ (I\otimes \tau)\circ m^*(f)\\
%&= i^*\circ (I\otimes \sigma)\circ m^*\circ i^*(\{g\mapsto \tau(\{h\mapsto f(gh)\})\}) \\
%&= i^*\circ (I\otimes \sigma)(\{(g,k)\mapsto \tau(\{h\mapsto f(k^{-1}g^{-1}h)\})\})\\
%&=\{g\mapsto \sigma(\{k\mapsto \tau(\{h\mapsto f(k^{-1}gh)\})\})\}\\
%&=\{g\mapsto(-1)^{|\tau||\sigma|}\tau(\{k\mapsto \sigma(\{h\mapsto f(h^{-1}gk)\})\})\},
%\end{align*} which is what we wanted to prove.
%%%%%%%%%%% BIS HIER
\end{proof}

\begin{Prop}\label{FormulaNablaBracket} Let $G$ be a Lie supergroup, equipped with a left-invariant graded Riemannian metric such that the induced scalar superproduct on $\mathfrak g$ is $\ad_{\mathfrak{g}}$-invariant. Then \begin{equation} \nabla_XY=\frac{1}{2}[X,Y] \quad \text{ and } \quad R(X,Y)Z=-\frac{1}{4}[[X,Y],Z]\end{equation} for all left-invariant vector fields $X,Y$.
\end{Prop}
\begin{Rem} Note that we do not yet prove the statement for bi-invariant metrics! 
\end{Rem}
\begin{proof} This is -- again apart from the additional signs -- the standard proof via (\ref{LeviCivita_Formula}). 
%%%%%%%%%%%% IM ARTIKEL WEG
%For left-invariant vector fields, (\ref{LeviCivita_Formula}) reduces to \begin{align} 2\left<\nabla_XY,Z\right>&= \left<[X,Y],Z\right>-(-1)^{|X|(|Y|+|Z|)}\left<[Y,Z],X\right> \nonumber\\
%&+ (-1)^{|Z|(|X|+|Y|)}\left<[Z,X],Y\right>.
%\end{align} 
%Using the $\ad_{\mathfrak{g}}$-invariance, we obtain
%\begin{align*}2\left<\nabla_XY,Z\right>&=\left<[X,Y],Z\right>-(-1)^{|X|(|Y|+|Z|)}\left<[Y,Z],X\right> -(-1)^{|Z||Y|}\left<X,[Z,Y]\right> \\
%&=\left<[X,Y],Z\right>,
%\end{align*} whereupon the first equation follows. The second one now directly follows from the first using the graded Jacobi identity (\ref{Jacobi}).
%%%%%%%%%%%% IM ARTIKEL WEG
%\begin{align*} R(X,Y)Z&=\nabla_X\nabla_YZ-(-1)^{|X||Y|}\nabla_Y\nabla_XZ-\nabla_{[X,Y]}Z \\
%&=\frac{1}{4}[X,[Y,Z]]-(-1)^{|X||Y|}\frac{1}{4}[Y,[X,Z]]-\frac{1}{2}[[X,Y],Z]\\
%&=-\frac{1}{4}[[X,Y],Z]
%\end{align*}by the graded Jacobi identity (\ref{Jacobi}). 
%%%%%%%%%%%% BIS HIER
\end{proof}

\begin{Lemma}\label{RightInvSpDiffAt_e} Let $G$ be a Lie supergroup with a left-invariant graded Riemannian metric $\left<\cdot,\cdot\right>$ such that the induced scalar superproduct on $\mathfrak g$ is $\ad_{\mathfrak{g}}$-invariant. Then $$\tau\left<X,Y\right>=0$$ for all right-invariant vector fields $X,Y$ and all tangent vectors $\tau\in T_eG$.
\end{Lemma}
\begin{proof} Let $\tau,\sigma$ and $\rho\in T_eG$ and denote by $X_\tau,X_\sigma,X_\rho$ and $Y_\tau,Y_\sigma,Y_\rho$ the corresponding left- and right-invariant vector fields. Then we have 
\begin{align*} &\tau\left<Y_\sigma,Y_\rho\right>=\left<\nabla_{X_\tau} Y_\sigma,Y_\rho\right>(e)+(-1)^{|\tau||\sigma|}\left<Y_\sigma,\nabla_{X_\tau}Y_\rho\right>(e)\\
&\quad=\left<(-1)^{|\tau||\sigma|}\nabla_{Y_\sigma}X_\tau+[X_\tau,Y_\sigma],X_\rho\right>(e)\\
&\quad\quad\quad\quad\quad\quad\quad\quad\quad\quad\quad\quad+(-1)^{|\tau||\sigma|}\left<X_\sigma,(-1)^{|\tau||\rho|}\nabla_{Y_\rho}X_\tau+[X_\tau,Y_\rho]\right>(e)\\
&\quad=(-1)^{|\tau||\sigma|}\left<\nabla_{X_\sigma}X_\tau,X_\rho\right>(e)+(-1)^{|\tau|(|\sigma|+|\rho|)}\left<X_\sigma,\nabla_{X_\rho}X_\tau\right>(e)\\
%&=(-1)^{|\tau||\sigma|}\frac{1}{2}\left(\left<[X_\sigma,X_\tau](e),X_\rho(e)\right>\right.\\
&\quad=-\frac{1}{2}\left(\left<[X_\tau,X_\sigma],X_\rho\right>(e)+(-1)^{|\tau||\sigma|}\left<X_\sigma,[X_\tau,X_\rho]\right>(e)\right)=0.
\end{align*} In the first and second line we used the properties of the Levi-Civita connection. Then we applied Lemma \ref{BracketLeftRight} and used the fact that for vector fields $Z$ and $W$, the value of $\nabla_Z W$ at some point depends only on the value of $Z$ at that point. In the last line we finally used Proposition \ref{FormulaNablaBracket} and the $\ad_{\mathfrak{g}}$-invariance of the induced scalar superproduct. 
\end{proof}
%\begin{proof} For $\tau\in T_eG$ let $X_\tau$, resp. $Y_\tau$ denote the unique left-, resp. right-invariant vector field whose value at $e$ is $\tau$. Since we know $[X_\tau,X_\sigma](e)=-[Y_\tau,Y_\sigma](e)$ for all $\tau,\sigma\in T_eG$ (cf. \cite{DelMor}, p.74), we conclude 
%\begin{align} \left<[Y_\tau,Y_\sigma],Y_\rho\right>(e)=-(-1)^{|\tau||\sigma|}\left<Y_\sigma,[Y_\tau,Y_\rho]\right> (e) \label{ad_inv_right}
%\end{align} for all tangent vectors $\tau,\sigma$ and $\rho$.
%The formula (\ref{LeviCivita_Formula}), together with (\ref{ad_inv_right}), yields
%\begin{align*} 2\left<\nabla_{Y_\tau} {Y_\sigma},{Y_\rho}\right>(e) &= {Y_\tau}\left<{Y_\sigma},{Y_\rho}\right>(e) -(-1)^{|\rho|(|\tau|+|\sigma|)}{Y_\rho}\left< {Y_\tau},{Y_\sigma}\right>(e) 
%\nonumber \\ 
% &+(-1)^{|\tau|(|\sigma|+|\rho|)}{Y_\sigma}\left< {Y_\rho},{Y_\tau}\right>(e)+\left< %[{Y_\tau},{Y_\sigma}],{Y_\rho}\right>(e).
%\end{align*}
%\marginpar{nicht fertig}
%\end{proof}

\begin{Thm} \label{BiAd} Let $G$ be a Lie supergroup with a left-invariant graded Riemannian metric $\left<\cdot,\cdot\right>$. Then the metric is bi-invariant if and only if $\left<\cdot,\cdot\right>_e$ is $\Ad_G$-invariant.
\end{Thm}
%%%%%%%%%%%%%% IM ARTIKEL WEG
%\begin{Rem} We emphasize that in the case of $G$ being connected this means that $\left<\cdot,\cdot\right>$ is bi-invariant if and only if $\left<\cdot,\cdot\right>_e$ is $\ad_{\mathfrak{g}}$--invariant.
%\end{Rem}
\begin{proof}
If the metric is left-invariant (resp. right-invariant), left (resp. right) trans\-lations with elements of $G_{\red}$ are isometries of $M$; 
%%%%%%%%%%%% IM ARTIKEL WEG
%\begin{align*}l_g^*\left<l_{g^{-1}}^*\circ X\circ l_g^*,l_{g^{-1}}^*\circ Y\circ l_g^*\right>&\overset{(\ref{left-invariant_translation})}{=}l_g^*\left<l_{g^{-1}}^*\circ l_g^*\circ X,l_{g^{-1}}\circ l_g^*\circ Y\right>\\
%&=l_g^*\left<X,Y\right>=\left<X,Y\right>\end{align*} 
%for all left-invariant vector fields $X,Y$ on $G$ (or analogously for right-invariant vector fields). 
%%%%%%%%%%%% BIS HIER
thus, if $\left<\cdot,\cdot\right>$ is bi-invariant, $\left<\cdot,\cdot\right>_e$ is $\Ad_{G_{\red}}$-invariant. To prove the $\ad_{\mathfrak{g}}$-invariance, calculate for \mbox{$\tau,\sigma,\rho\in T_eG$} as follows ($X_\tau,X_\sigma,X_\rho$ and $Y_\tau,Y_\sigma,Y_\rho$ are again the corresponding left- and right-invariant vector fields):
\begin{align*}
&\left<[X_\tau,X_\sigma],X_\rho\right>_e=\left<\nabla_{X_\tau} X_\sigma-(-1)^{|\tau||\sigma|}\nabla_{X_\sigma}X_\tau,X_\rho\right>_e\\
&\quad=\left(X_\tau\left<X_\sigma,X_\rho\right>(e)-(-1)^{|\tau||\sigma|}\left<X_\sigma,\nabla_{X_\tau}X_\rho\right>_e\right)-(-1)^{|\tau||\sigma|}\left<\nabla_{Y_\sigma}X_\tau, X_\rho\right>(e)\\
&\quad=-(-1)^{|\tau||\sigma|}\left<X_\sigma,\nabla_{X_\tau}X_\rho\right>_e-\left<\nabla_{X_\tau}Y_\sigma+(-1)^{|\tau||\sigma|}[Y_\sigma,X_\tau],Y_\rho\right>(e)\\
&\quad=-(-1)^{|\tau||\sigma|}\left<X_\sigma,\nabla_{X_\tau}X_\rho\right>_e-\left(X_\tau\left<Y_\sigma,Y_\rho\right>(e)-(-1)^{|\tau||\sigma|}\left<Y_\sigma,\nabla_{X_\tau}Y_\rho\right>(e)\right)\\
&\quad=-(-1)^{|\tau||\sigma|}\left<X_\sigma,\nabla_{X_\tau}X_\rho\right>_e+(-1)^{|\tau|(|\sigma|+|\rho|)}\left<X_\sigma,\nabla_{Y_\rho}X_\tau\right>(e)\\
&\quad=-(-1)^{|\tau||\sigma|}\left<X_\sigma,\nabla_{X_\tau}X_\rho-(-1)^{|\tau||\rho|}\nabla_{X_\rho}X_\tau\right>_e\\
&\quad=-(-1)^{|\tau||\sigma|}\left<X_\sigma,[X_\tau,X_\rho]\right>_e.
\end{align*}
Note that we used Lemma \ref{BracketLeftRight} for the fourth and the right-invariance of $\left<\cdot,\cdot\right>$ for the fifth equality. 

Let now $\left<\cdot,\cdot\right>_e$ be $\Ad_G$-invariant. We want to show the bi-invariance of $\left<\cdot,\cdot\right>$. The left-invariance of $\left<\cdot,\cdot\right>$, together with the $\Ad_{G_{\red}}$-invariance of $\left<\cdot,\cdot\right>_e$ shows that for all right-invariant vector fields $X$ and $Y$, the $\shC^\infty$-part of the superfunction $\left<X,Y\right>$ is constant. We have to show that the superfunction itself is constant.

Fix a basis $\{Z_\eta\}=\{Z_i,Z_\alpha\}$ consisting of homogeneous right-invariant vector fields; assume that the value of $Z_\eta$ at $e$ is $\partial_\eta$. We write $$Z_\eta=\sum_iz_\eta^i\partial_i+\sum_\alpha z_\eta^\alpha\partial_\alpha$$ for superfunctions $z_\eta^i$ and $z_\eta^\alpha$. 
First of all, let us have a look at the case that $X$ is even and $Y$ is odd. We write $$\left<X,Y\right>=\sum_{\alpha}f_\alpha\xi_\alpha + O(\xi^3)$$ for some $\shC^\infty$-functions $f_\alpha$. Our first goal is to show that $f_\alpha=0$. On the one hand $$Z_\beta\left<X,Y\right>=\sum_{\alpha,i}z_\beta^i(\partial_i f_\alpha)\xi_\alpha+\sum_\alpha z_\beta^\alpha f_\alpha+O(\xi^2)=\sum_\alpha \widetilde{z_\beta^\alpha}f_\alpha+O(\xi^2).$$ On the other hand, $$Z_\beta\left<X,Y\right>=\left<\nabla_{Z_\beta}X,Y\right>+\left<X,\nabla_{Z_\beta}Y\right>$$ is the sum of a constant function and nilpotent terms. Consequently, $$\sum_\alpha \widetilde{z_\beta^\alpha}f_\alpha\equiv Z_\beta\left<X,Y\right>(e)=0$$ according to Lemma \ref{RightInvSpDiffAt_e}. Since the matrix $(\widetilde{z_\beta^\alpha})_{\beta,\alpha}$ is invertible (as $\{Z_\eta\}$ is a basis) we conclude that all the $f_\alpha$ must vanish.

Let us now have a look at the case when $X$ and $Y$ are both even or both odd. Then %$$\left<X,Y\right>=\text{const}+\sum_{\alpha<\beta}f_{\alpha\beta}\xi_\alpha\xi_\beta+O(\xi^4)$$ for some $\shC^\infty$-functions $f_{\alpha\beta}$.
$$\left<X,Y\right>=\text{const}+f_2+O(\xi^4)$$ for some homogeneous function $f_2$ of degree $2$ (in the chosen coordinates). Then for any $\gamma$,
$$Z_\gamma\left<X,Y\right>=\sum_\alpha \widetilde{z_\gamma^\alpha} (\partial_\alpha f_2) +O(\xi^3).$$ But because $$Z_\gamma\left<X,Y\right>=\left<\nabla_{Z_\gamma}X,Y\right>+(-1)^{|X|}\left<X,\nabla_{Z_\gamma}Y\right>$$ consists only of terms of degree at least three by what we have shown before, and because the matrix $(\widetilde{z_\gamma^\alpha})$ is invertible, we see that $\partial_\alpha f_2=0$ for any $\alpha$. The superfunction $f_2$ is homogeneous of degree $2$, so $f_2=0$.

We proceed by induction: Consider the case of $X$ even and $Y$ odd, write %$$\left<X,Y\right>=\sum_{\alpha_1<\ldots<\alpha_{2k+1}}f_{\alpha_1\ldots\alpha_{2k+1}}\xi_{\alpha_1}\ldots\xi_{\alpha_{2k+1}}+O(\xi^{2k+3}),$$ 
$$\left<X,Y\right>=f_{2k+1} + O(\xi^{2k+3})$$ for some homogeneous superfunction $f_{2k+1}$ of degree $2k+1$, apply all the $Z_\beta$ and conclude that $f_{2k+1}=0$. The analogous argument as above for the case of $X$ and $Y$ of the same parity finishes the induction.
%\marginpar{und?}
%For the converse we have to take a closer look on how the adjoint representation, together with the Lie structure on $G_{\red}$, determine the supergroup structure on $G$. \marginpar{?????? proof}
\end{proof}

\begin{Cor}\label{FormulaNablaR} Let $G$ be a Lie supergroup with a bi-invariant graded Riemannian metric. Then \begin{equation} \nabla_XY=\frac{1}{2}[X,Y] \quad \text{ and } \quad R(X,Y)Z=-\frac{1}{4}[[X,Y],Z]\end{equation} for all left-invariant vector fields $X,Y$.
\end{Cor}

As the last result in this section, we state the following easy generalization of Proposition 1.6 b) of \cite{CahPar} -- it will reduce the amount of calculation needed for the last example in \ref{Section_Examples}. Up to signs, the proof is the same as there. This result also appears as Theorem 4 in \cite{Cortes}; note that there the assumption of the adjoint representation being faithful is missing.
\begin{Prop}\label{ScalarSP_OddPart} Let ${\mathfrak g}={\mathfrak g}_0\oplus {\mathfrak g}_1$ be a Lie superalgebra such that $[{\mathfrak g}_1,{\mathfrak g}_1]={\mathfrak g}_0$ and the adjoint representation of ${\mathfrak g}_0$ on ${\mathfrak g}_1$ is faithful. Then any non-degenerate $\ad_{{\mathfrak g}_0}$-invariant skew-symmetric bilinear form $\left<\cdot,\cdot\right>_1$ on ${\mathfrak g}_1$ uniquely extends to an  $\ad_{{\mathfrak g}}$-invariant scalar superproduct $\left<\cdot,\cdot\right>$ on ${\mathfrak g}$.
\end{Prop}

\subsection{The Killing Form} \label{SectionKilling}

%Most part of this section is taken from \cite{Kac} and \cite{Scheunert}. 
%Let $K$ denote either the field $\R$ or $\C$ and 
Let ${\mathfrak g}$ be a finite-dimensional Lie superalgebra. We define the {\it Killing form} of ${\mathfrak g}$ to be \begin{equation}B(X,Y):=\str(\ad_X\circ \ad_Y),\end{equation} where $\str$ is the  supertrace. The Killing form is an even graded-symmetric $\ad_{\mathfrak{g}}$-invariant bilinear form; the $\ad_{\mathfrak{g}}$-invariance meaning \begin{equation}B([X,Y],Z)+(-1)^{|X||Y|}B(Y,[X,Z])=0\end{equation} for all $X,Y,Z\in {\mathfrak g}$.% If $G$ is a connected Lie supergroup with Lie algebra ${\mathfrak g}$, the Killing form is $\ad$-invariant.

It is known (see e.g. \cite{Scheunert}), that on ${\mathfrak {sl}}(n|m)$, the Killing form is given by $$B(X,Y)=2(n-m)\,\str(XY);$$ it is non-degenerate for $n\neq m$ but identically zero for $n=m$. Nevertheless, the supertrace $(X,Y)\mapsto \str(XY)$ induces a non-degenerate even graded-symmetric $\ad_{\mathfrak{sl}(n|m)}$-invariant bilinear form on ${\mathfrak{sl}}(n|n)/{K\cdot I_{2n}},$ where $I_{2n}$ is the $2n\times 2n$ unit matrix.

On ${\mathfrak {osp}}(n|2m)$, the Killing form is $$B(X,Y)=(n-2m-2)\,\str(XY);$$ it is non-degenerate if $n\neq 2m+2$ and identically zero for $n=2m+2$. Here, the supertrace is non-degenerate and $\ad_{\osp(n|2m)}$-invariant for all $n$ and $m$.

%If $K$ is algebraically closed the Lie superalgebras with non-degenerate Killing form are classified (\cite{Scheunert}, corollaries to Th. 1, Ch. II, $\S$3 and Th. 1, Ch. II, $\S$5): The Killing form of a Lie superalgebra ${\mathfrak g}$ over $K$ is non-degenerate if and only if ${\mathfrak g}$ is the direct product of classical simple Lie superalgebras whose Killing forms are non-degenerate. Those are exactly the simple Lie algebras, the Lie superalgebras \begin{align*}
%&{\mathfrak{sl}}(n|m), \quad n,m\ge 1, n\neq m, \\
%&{\mathfrak{osp}}(n|2m),\quad n,m\ge 1,n\neq 2m+2, 
%\end{align*} and the exceptional Lie superalgebras ${\mathfrak f}_4$ and ${\mathfrak g}_3$ (see \cite{Scheunert} for the definition of those).
 
Note that there exist Lie superalgebras with degenerate (even vanishing) Killing form that nevertheless admit an $\ad$-invariant scalar superproduct. 

\subsection{Homogeneous Superspaces}

If $G$ is a Lie supergroup and $H$ a closed Lie subsupergroup (a Lie subsupergroup such that $G_{\red}$ is closed in $H_{\red}$), we know from the standard theory that $G_{\red}/H_{\red}$ is a manifold. Consider the canonical projections $\pi:G_{\red} \to G_{\red}/H_{\red}$ and ${\rm pr}_1:G\times H\to G$ and the right action of $H$ on $G$, $\Phi=(\phi,\phi^*):G\times H\to G$, i.e.~the composition of the multiplication morphism with the inclusion of $H$ into $G$. We equip $G_{\red}/H_{\red}$ with the sheaf of $H$-invariant superfunctions:
\begin{equation}\shO_{G/H}(U):=\{f\in \shO_G(\pi^{-1}(U))\mid \phi^*f={\rm pr}_1^*f\}.\label{SheafOfHomSpace}\end{equation} Then $$G/H:=(G_{\red}/H_{\red},\shO_{G/H})$$ is a supermanifold, see \cite{Kostant}, Theorem 3.9. There is a canonical morphism of supermanifolds $G\to G/H$, where the sheaf morphism is given by inclusion. %For a more general discussion of supermanifolds of orbits, see \cite{Stavracou}.
The surjective morphism ${\mathfrak g}\cong T_eG\to T_{H}G/H$ has kernel $\mathfrak h$ and thus yields a canonical identification $$T_{H}G/H\cong {\mathfrak g}/{\mathfrak h}.$$ The adjoint representation of $G$ restricts to a well-defined representation $\Ad_H$ of $H$ on ${\mathfrak g}/{\mathfrak h}$, i.e.~compatible representations $$\Ad_{H_{\red}}:H_{\red}\times {\mathfrak g}/{\mathfrak h}\text{ and } \ad_{\mathfrak h}:{\mathfrak h}\times {\mathfrak g}/{\mathfrak h}\to {\mathfrak g}/{\mathfrak h}.$$ 

If $\rho:G\times M\to M$ is an action of a Lie supergroup $G$ on a supermanifold $M$ and $p\in M_{\red}$, we define the {\it isotropy group} $G_p$ %via the functor of points: $G_x(S):=G(S)_{x_S}$, the stabilizer of the action of $G(S)$ on $M(S)$ of the constant mapping $x_S$. This functor can be shown to be representable \marginpar{really?}and thus defines a Lie subsupergroup of $G$. The reduced Lie group of the stabilizer is the stabilizer of the reduced action: $$(G_x)_{\red}=(G_{\red})_x.$$
to be the Lie subsupergroup of $G$ given by the Harish-Chandra pair
$$((G_{\red})_p,{\mathfrak g}_p),$$ where \begin{equation}{\mathfrak g}_p:=\{X\in {\mathfrak g}\mid ev_p\circ (X_e\otimes I)\circ \rho^*=0\}.\label{IsotropyAlgebra}\end{equation} %\marginpar{von Onishchik abgeschrieben (S. Lie memorial...)} 
Here, $ev_p:\shT_M(M)\to T_pM$ is evaluation at $p$; cf. also \ref{Section_Action}. In other words, ${\mathfrak g}_p$ consists of those left-invariant vector fields  $X$ on $G$ such that the infinitesimal action of $X$ at $p$ is trivial.

As usual we have the {\it isotropy representation} of $G_p$ on $T_pM$: On the level of the underlying Lie group, $g\in (G_{\red})_p$ acts on $T_pM$ in an obvious manner; the Lie superalgebra ${\mathfrak g}_p$ acts on %\marginpar{sign correct?}
$T_pM$ via \begin{equation}X\cdot v=-[(X_e\otimes I)\circ \rho^*,v].\label{IsotropyRepresentation}\end{equation} Note that this is a well-defined action compatible with the $(G_{\red})_p$-action, so they fit together to a representation of the Lie supergroup $G_p$.

Let $\rho:G\times M\to M$ be an action on a supermanifold $M$ and fix $p\in M_{\red}$. We can define the {\it orbit} $G\cdot p$ as follows: The usual orbit $G_{\red}\cdot p\subset M_{\red}$ inherits its differentiable structure via the canonical mapping $j:G_{\red}/(G_{\red})_p\to G_{\red}\cdot p$. The supermanifold $G\cdot p$ is now defined to be $$G\cdot p:=(G_{\red}\cdot p,j_*\shO_{G/G_p});$$ it is no surprise that also the superstructure has to be passed over from $G/G_p$. If we denote the inclusion $G_{\red}\cdot p\to M_{\red}$ by $i$, the orbit map $$\rho_p:G\overset{\id\times p}{\longrightarrow} G\times M\overset{\rho}{\longrightarrow} M$$ of the point $p$ gives a well-defined sheaf morphism $$\shO_M\to i_*\shO_{G\cdot p}=(i\circ j)_*\shO_{G/G_p}.$$ Note that the orbits $G\cdot p$ become submanifolds of $M$ via this inclusion morphism \mbox{$i=(i,\rho_p^*):G\cdot p\to M$}.

An action $\rho:G\times M\to M$ is said to be {\it transitive} (cf. \cite{OnishchikFlag}, p. 295) if the reduced action $\tilde{\rho}:G_{\red}\times M_{\red}\to M_{\red}$ is transitive and if for all $p\in M_{\red}$, the mapping \begin{equation}{\mathfrak g}\to T_pM;\quad X\mapsto ev_p\circ (X_e\otimes I)\circ \rho^*=X_e\circ \rho_p^*\label{IsomorphismGTpM}\end{equation} is surjective. It is known from the standard theory that in the case of a connected supermanifold, the first condition is equivalent to the even part of the second condition, i.e.~the surjectivity of the mapping ${\mathfrak g}_0\to (T_pM)_0$. We say that a supermanifold is {\it $G$-homogeneous} if it is acted on transitively by $G$; we say it is {\it homogeneous} if it is $G$-homogeneous for some Lie supergroup $G$. %\marginpar{cf. Onishchik, Homogeneous supermanifolds...}
%The mapping (\ref{IsomorphismGTxM}) descends to an isomorphism ${\mathfrak g}/{\mathfrak g}_x\to T_xM$.
If $M$ is already assumed to be Riemannian, it is {\it homogeneous} if it is $I(M)$-homogeneous.

Note that if $M$ is $G$-homogeneous, the natural morphism $G\cdot p\to M$ is an isomorphism.%: On the level of the underlying action it is clear, and the rest follows from the definition of $G$-homogeneity and the definition of ${\mathfrak g}_p$.

\begin{Prop} \label{IsotropyAdjointEquiv} Let $M$ be $G$-homogeneous. Then the isotropy representation of $G_p$ on $T_p M$ is equivalent to the adjoint representation on ${\mathfrak g}/{\mathfrak g}_p$ via the natural isomorphism ${\mathfrak g}/{\mathfrak g}_p\to T_pM$.
\end{Prop}
\begin{proof} First we have a look at the representations of the underlying Lie group. For any $g\in (G_{\red})_p$ and any $X\in {\mathfrak g}$, we have $$(\Ad_g X)_e\circ \rho_p^*=X_e\circ l_g^*\circ r_{g^{-1}}^*\circ \rho_p^*=X_e\circ l_g^*\circ \rho_p^*= X_e\circ \rho_p^*\circ g^*= dg(X_e\circ \rho_p^*),$$ where $l_g$ and $r_g:G\to G$ are left and right multiplication with $g$, and the diffeomorphism of $M$ induced by the action of $g$ is again denoted by $g:M\to M$. %, which in symbolic notation works as follows: for a superfunction $f$ on $M$, we have
%\begin{align*} ev_p\circ i^*\circ(\Ad_g X\otimes I)\circ \phi^* f&=\Ad_g X(\{h\mapsto f(hp)\})(e)\\
%&=X(\{h\mapsto f(ghg^{-1}p)\})(e)\\
%&=X(\{h\mapsto f(ghp)\})(e)\\
%&=ev_p\circ i^*\circ (X\otimes I)\circ \phi^*\circ g^* f.
%\end{align*} 

The equivalence of the representations at the level of Lie superalgebras is the equation
$$[Y,X]_e\circ \rho_p^*=-[(Y_e\otimes I)\circ \rho^*, X_e\circ \rho_p^*]$$ for all $X\in {\mathfrak g}$ and all $Y\in {\mathfrak g}_p$, cf. (\ref{IsotropyRepresentation}), which follows directly from (\ref{BracketsSign}).
\end{proof}

\subsection{Invariant Metrics on Homogeneous Superspaces}

In the standard theory a Riemannian metric on a $G$-homogeneous space $M$ is called $G$-invariant if every $g\in G$ acts on $M$ by isometries. The characterization of actions in the world of supermanifolds in terms of Harish-Chandra pairs immediately motivates the following definition:

A graded Riemannian metric on a $G$-homogeneous superspace $M$ is {\it $G$-invariant} if every $g\in G_{\red}$ acts on $M$ by isometries and the image of the morphism ${\mathfrak g}\to \shT_M(M)^\circ$  lies in the subalgebra of graded Killing fields. 

\begin{Thm}\label{InvMetricHomSpace} Let $M$ be a $G$-homogeneous superspace, fix $p\in M_{\red}$ and let $H=G_p$. Then there is a 1-1-correspondence between $G$-invariant graded Riemannian metrics on $M$ and $\Ad_H$-invariant scalar superproducts on the space \mbox{${\mathfrak g}/{\mathfrak h}\simeq T_{H} G/H\simeq T_pM$}.
\end{Thm}
\begin{Rem} By definition, $\Ad_H$-invariance means $\Ad_{H_{\red}}$- and $\ad_{\mathfrak{h}}$-invariance. If $H$ is connected, $\ad_{\mathfrak h}$-invariance implies $\Ad_{H_{\red}}$-invariance. Nevertheless, the converse is not true -- $\Ad_{H_{\red}}$-invariance only implies $\ad_{{\mathfrak h}_0}$-invariance.
\end{Rem}
\begin{proof} First of all, take a $G$-invariant graded Riemannian metric $\left<\cdot,\cdot\right>$ on $M$; we have to show that the induced scalar superproduct on ${\mathfrak g}/{\mathfrak h}$ is $\Ad_{H_{\red}}$- and $\ad_{\mathfrak h}$-invariant.

Since the action of any element $g\in G_{\red}$ and thus in particular of any element $h\in H_{\red}$ is isometric, Proposition \ref{IsotropyAdjointEquiv} immediately shows that $\left<\cdot,\cdot\right>_H$ is \mbox{$\Ad_{H_{\red}}$-in}\-variant.

Let $X\in {\mathfrak h}$. If the action is denoted by $\rho$, then $\bar{X}:=(X_e\otimes I)\circ \rho^*$ is a Killing field with value $0$ at $p$. Pick $v,w\in T_pM$ and extend them to vector fields $Y,Z$ on a neighbourhood of $p$. Then
\begin{align*}0=\bar{X}\left<Y,Z\right>(p)&=\left<[\bar{X},Y],Z\right>(p)+(-1)^{|\bar{X}||Y|}\left<Y,[\bar{X},Z]\right>(p)\\
&=\left<[\bar{X},v],w\right>_p+(-1)^{|\bar{X}||v|}\left<v,[\bar{X},w]\right>_p;\end{align*} note that the expressions in the second line are well-defined since $\bar{X}$ vanishes at $p$. Translating the situation to ${\mathfrak g}/{\mathfrak h}$ with the help of Proposition \ref{IsotropyAdjointEquiv} gives the $\ad_{\mathfrak h}$-invariance.

We now show the injectivity of the correspondence: Assume that there are $G$-invariant graded Riemannian metrics $\left<\cdot,\cdot\right>_1$ and $\left<\cdot,\cdot\right>_2$ that induce the same scalar superproduct on $T_pM$. The fact that every $g\in G_{\red}$ acts on $M$ by isometries with respect to both metrics immediately implies that they induce the same scalar superproduct on every tangent space. But that does not suffice to show that the metrics are equal. To achieve that, take two local vector fields $Y,Z$ on $M$. We know that $$\left<Y,Z\right>_1(q)=\left<Y,Z\right>_2(q)$$ for any $q$. Furthermore, for any Killing field $X_1$, \begin{align*}X_1\left<Y,Z\right>_1(q)&=\left<[X_1,Y],Z\right>_1(q)+(-1)^{|X_1||Y|}\left<Y,[X_1,Z]\right>_1(q)\\
&=\left<[X_1,Y],Z\right>_2(q)+(-1)^{|X_1||Y|}\left<Y,[X_1,Z]\right>_2(q)=X_1\left<Y,Z\right>_2(q).\end{align*} Applying further Killing fields $X_2,\ldots,X_n$, we see analogously
$$X_n\ldots X_1\left<Y,Z\right>_1(q)=X_n\ldots X_1\left<Y,Z\right>_2(q),$$ since there only appear scalar products of some iterated Lie brackets, evaluated at $q$. The $G$-homogeneity of $M$ now gives the existence of enough Killing fields to conclude $\left<Y,Z\right>_1=\left<Y,Z\right>_2$.

It remains to show the surjectivity of the correspondence. Let $\left<\cdot,\cdot\right>_H$ be an $\Ad_H$-invariant scalar superproduct on ${\mathfrak g}/{\mathfrak h}$, which by Proposition \ref{IsotropyAdjointEquiv} is the same as a scalar superproduct on $T_pM$, invariant under the isotropy representation of $G_p$, i.e.~invariant under $(G_{\red})_p$ and ${\mathfrak g}_p$, where the second condition is $\left<X\cdot v,w\right>_p+(-1)^{|X||v|}\left<v,X\cdot w\right>_p=0$. The action of $G_{\red}$ now induces well-defined scalar superproducts $\left<\cdot,\cdot\right>_q$ on $T_qM$ for all $q\in M_{\red}$: for $g\in G_{\red}$ with $gq=p$ we set $\left<v,w\right>_q:=\left<dg(v),dg(w)\right>_p$. These are invariant unter the action of the respective isotropy groups -- on the level of Lie groups this follows directly from $(G_{\red})_q=g^{-1}\cdot(G_{\red})_p\cdot g$; on the level of Lie superalgebras we calculate for $X\in {\mathfrak g}_q$ as follows:  
\begin{align*}\left<X\cdot v,w\right>_q&=-\left<[\bar{X},v],w\right>_q=-\left<dg[\bar{X},v],dg(w)\right>_p\\
&=-\left<[dg(\bar{X}),dg(v)],dg(w)\right>_p=(-1)^{|X||v|}\left<dg(v),[dg(\bar{X}),dg(w)]\right>_p\\
&=-(-1)^{|X||v|}\left<v,X\cdot w\right>_q.\end{align*}

Let $X_1,\ldots,X_n$ be a basis of the Lie algebra ${\mathfrak g}$ and consider the associated vector fields $\bar{X}_i=((X_i)_e\otimes I)\circ \rho^*$. If $(\eta_i)$ are coordinates on some open set $U$, we can express $\partial_j$ as a linear combination of the $\bar{X}_i$, since $M$ is $G$-homogeneous: \begin{equation}\partial_j=\sum_i f^j_i \bar{X}_i.\label{KoordInKilling}\end{equation} We may assume that \begin{equation} |\partial_j|=|f^j_i|+|\bar{X}_i| \label{invMetric_parity_convention}
\end{equation} for all $i$ and $j$. Of course, there is no unique way of doing so, but if \begin{equation}\sum_i f_i \bar{X}_i=\sum_i g_i \bar{X}_i\label{TwoReprKilling}\end{equation} are two representations of the same vector field, we have $$\sum_i (f_i(q)-g_i(q))X_i\in {\mathfrak g}_q$$ for all $q\in U$ by definition of ${\mathfrak g}_q$ (\ref{IsotropyAlgebra}). 

Distinguish now between even and odd coordinates $x_i$ and $\xi_\alpha$ and take two vector fields $Y$ and $Z$ on $U$. Our strategy is as follows: Under the assumption of its existence we compute $\left<Y,Z\right>$ in terms of the scalar superproducts at all points $q\in U$; then we show that the resulting expressions provide us with a well-defined $G$-invariant metric. \begin{equation}\left<Y,Z\right>=\left<Y,Z\right>(q)+\sum_\alpha (\partial_\alpha \left<Y,Z\right>)(q)\xi_\alpha+\sum_{\alpha<\beta}(\partial_\alpha\partial_\beta \left<Y,Z\right>)(q)\xi_\beta\xi_\alpha+\ldots\label{MetricCoordinates}\end{equation}
Choosing some representation of $\partial_\alpha$ in terms of the $\bar{X}_i$, (\ref{KoordInKilling}), we calculate
\begin{align*}(\partial_\alpha \left<Y,Z\right>)(q)%&=\sum_i f^\alpha_i(q) (\bar{X}_i\left<Y,Z\right>)(q)\\
&=\sum_i f^\alpha_i(q)\left(\left<[\bar{X}_i,Y],Z\right>(q)+(-1)^{|\bar{X}_i||Y|}\left<Y,[\bar{X}_i,Z]\right>(q)\right).\end{align*}
The right-hand side now defines a $\shC^\infty$-function $G^{Y,Z}_\alpha(q)$; it is independent of the representation of $\partial_\alpha$ in the $\bar{X}_i$ since the scalar superproduct at $q$ is $G_q$-invariant.%\marginpar{ausf\"uhren} 

Again assuming the existence of $\left<Y,Z\right>$, we compute the part of degree two:
\begin{align*} (\partial_\alpha &\partial_\beta\left<Y,Z\right>)(q)\\
&=\partial_\alpha\left(\sum_i f_i^\beta \left(\left<[\bar{X}_i,Y],Z\right>+(-1)^{|\bar{X}_i||Y|}\left<Y,[\bar{X}_i,Z]\right>\right)\right)(q) \\
&=\sum_{i} (\partial_\alpha f^\beta_i)(q)\left(\left<[\bar{X}_i,Y],Z\right>(q)+(-1)^{|\bar{X}_i||Y|}\left<Y,[\bar{X}_i,Z]\right>(q)\right)\\
&\quad+\sum_{i,j} (-1)^{|f^\beta_i|}f^\beta_i(q) f^\alpha_j(q)\bigg(\left<[\bar{X}_j,[\bar{X}_i,Y]],Z\right>(q)\\
&\quad+(-1)^{|\bar{X}_j|(|\bar{X}_i|+|Y|)}\left<[\bar{X}_i,Y],[\bar{X}_j,Z]\right>(q)+(-1)^{|\bar{X}_i||Y|} \left<[\bar{X}_j,Y],[\bar{X}_i,Z]\right>(q)\\
&\quad+(-1)^{|Y|(|\bar{X}_i|+|\bar{X}_j|)}\left<Y,[\bar{X}_j,[\bar{X}_i,Z]]\right>(q)\bigg)
\end{align*} We define a $\shC^\infty$-function $G^{Y,Z}_{\alpha\beta}(q)$ by this expression and have to show that it is independent of the choice of representation (\ref{KoordInKilling}). Taking a second representation as in (\ref{TwoReprKilling}) we have 
\begin{align*}
0&=[\partial_\beta,\partial_\alpha]=\partial_\beta\circ\partial_\alpha +\partial_\alpha\circ \partial_\beta\\
&=\left(\sum_i f^\beta_i\bar{X}_i\right)\left(\sum_j g^\alpha_j\bar{X}_j\right)+\left(\sum_j g^\alpha_j\bar{X}_j\right)\left(\sum_i f^\beta_i\bar{X}_i\right)\\
&=\sum_j (\partial_\beta g^\alpha_j)\bar{X}_j+\sum_i (\partial_\alpha f^\beta_i) \bar{X}_i+\sum_{i,j} (-1)^{|\bar{X}_i||g^\alpha_j|}
f^\beta_ig^\alpha_j [\bar{X}_i,\bar{X}_j],\end{align*} where we used the parity convention (\ref{invMetric_parity_convention}). In other words,
\begin{equation} \sum_i \left((\partial_\beta g^\alpha_j)(q)+(\partial_\alpha f^\beta_i)(q)\right)X_i+\sum_{i,j}(-1)^{|X_i|}f^\beta_i(q) g^\alpha_j(q)[X_i,X_j]\in {\mathfrak g}_q
\end{equation} for all $q$. Using the $G_q$-invariance of the scalar superproduct at $q$, a short calculation shows that $G^{Y,Z}_{\alpha\beta}$ is independent of the choice of representation.

Similar but longer calculations give well-defined smooth functions $G^{Y,Z}_{\alpha_1\ldots\alpha_k}$ for any $k$-tuple $\alpha_1<\ldots<\alpha_k$. Then we can define a $G$-invariant metric by (\ref{MetricCoordinates}).%\marginpar{alles nachrechnen}
\end{proof}

\subsection{Riemannian Symmetric Superspaces}\label{Section_RSSS}

A conventional Riemannian manifold $M$ is called a Riemannian symmetric space if for every point $p\in M$ there exists an isometry $s_p$ of $M$ with $s_p(p)=p$ and $d_ps_p=-\id_{T_pM}$. Translating this into the world of supermanifolds we arrive at the following definition deviating from the usual one by the additional infinitesimal odd part: A Riemannian supermanifold $M$ is called {\it symmetric} or a {\it (Riemannian) symmetric superspace} if for any point $p$ there exists an iso\-metry $s_p$ of $M$ with $s_p(p)=p$ and $d_ps_p=-\id_{T_pM}$ and if for any odd tangent vector $\tau\in (T_pM)_1$ there exists a Killing vector field $S_\tau$ on $M$ with $(S_\tau)_p=\tau$ and $(\nabla S_\tau)(p)=0$. In the standard theory, the Killing vector fields $X$ with $(\nabla X)(p)=0$ are exactly those defined by transvections (which in turn are constructed via the geodesic symmetries) so the existence of Killing fields of this type is the correct infinitesimal counterpart of the existence of the geodesic symmetries.
%A Riemannian supermanifold $M$ is called {\it symmetric} or a {\it (Riemannian) symmetric superspace} if $M$ is homogeneous and there is a point $p\in M_{\red}$ with an isometry $\psi:M\to M$ such that $\psi(p)=p$ and $d_p\psi=-\Id_{T_pM}$. Note that then $\psi^2=\Id_M$ because of Proposition \ref{IsometryDetermined}. We impose homogeneity in the definition because of the difficulty of generalizing definitions containing conditions like "for all points..." occuring again. \marginpar{schoener: ohne Homog...} Clearly, if $M$ is homogeneous, such an isometry $\psi$ exists for any point $q\in M$.

\begin{Rem} Infinitesimal versions of this definition already exist in the mathema\-tical literature, see e.g.~\cite{Cortes} or \cite{Serganova}.
\end{Rem}

%A Riemannian supermanifold $M$ is called {\it locally symmetric} if the curvature tensor of $M$ is parallel (already defined in \cite{DeWitt}, Exercise 3.13 for the case of pseudo-Riemannian metrics).

Just like in the standard theory, a symmetric superspace is homogeneous (the surjectivity of the mappings $\{\text{odd Killing fields}\}\to (T_pM)_1$ is trivially fulfilled).

Let $M$ be a symmetric superspace and set $G=I_0(M)$, the identity component of the isometry group of $M$. Let $K=G_p$, the isotropy group of some point $p$. Conjugation with $s_p$ induces a morphism $\sigma:G\to G$ with $d\sigma=\Ad_{s_p}:{\mathfrak g}\to {\mathfrak g}$. % -- formally identify $\psi$ with the constant mapping $\psi:\R^{0|0}\to I(M)$ with image $\psi$ and define $\sigma=m\circ (\Id\times m)\circ (\psi\times \Id \times \psi):I(M)=\R^{0|0}\times I(M)\times \R^{0|0}\to I(M)$. 
Clearly, $\sigma$ is involutive.

We define the fixed point group $G^\sigma$ of $\sigma$ to be the Lie subsupergroup of $G$ given by the Harish-Chandra pair $(G^{\sigma}_{\red},{\mathfrak g}^{\sigma})$ with \mbox{$G^\sigma_{\red}=\{x\in G_{\red}\mid \sigma(x)=x\}$} and \mbox{${\mathfrak g}^{\sigma}=\{X\in {\mathfrak g}\mid d\sigma(X)=X\}$}. 
\begin{Lemma} Under these conditions, we have $G^\sigma_0\subset K\subset G^\sigma.$
\end{Lemma}
\begin{proof} The first inclusion can be verified on the level of Lie algebras: If $X$ is a Killing field on $M$ such that $d\sigma(X)=X$, i.e.~$ds_p(X)=X$, the value of $X$ at $p$ clearly has to vanish since $d_ps_p=-\id$.

Let $\phi\in K_{\red}$. Then $\phi$ and $\sigma(\phi)=s_p\circ\phi\circ s_p$ are isometries of $M$ sending $p$ to itself and having the same differential at $p$; because of Proposition \ref{IsometryDetermined}, they are equal, i.e.~$\phi\in G^{\sigma}_{\red}$. 

We have to show that ${\mathfrak k}$, the Lie algebra of $K$, is contained in ${\mathfrak g}^{\sigma}$, so let $X\in {\mathfrak g}$ be a Killing field that vanishes at $p$. Then $ds_p(X)$ is a Killing field also vanishing at $p$ and satisfying \begin{equation}(\nabla_Y X)(p)=-d_ps_p((\nabla_Y X)(p))=-(\nabla_{ds_p(Y)}ds_p(X))(p)=(\nabla_Y ds_p(X))(p)\label{KillingCalc}\end{equation} for all $Y$; Proposition \ref{KillingDeterminedByPoint} thus yields $X=ds_p(X)$.
\end{proof}

Since $d\sigma$ is an involutive automorphism, ${\mathfrak g}$ splits as the sum of its $(+1)$- and $(-1)$-eigenspace -- we can write any $X\in {\mathfrak g}$ as $X=\frac{1}{2}(X+d\sigma(X))+\frac{1}{2}(X-d\sigma(X))$. We showed that the $+1$-eigenspace coincides with the Lie algebra of $K$, so $${\mathfrak g}={\mathfrak k}\oplus{\mathfrak p}$$ with ${\mathfrak p}=\{X\in {\mathfrak g}\mid d\sigma(X)=-X\}$. The usual relations $[{\mathfrak k},{\mathfrak k}]\subset {\mathfrak k},\, [{\mathfrak k},{\mathfrak p}]\subset {\mathfrak p}$ and $[{\mathfrak p},{\mathfrak p}]\subset {\mathfrak k}$ hold.
\begin{Lemma} The space ${\mathfrak p}$ is the space of all Killing vector fields $X$ on $M$ such that $(\nabla X)(p)=0$.
\end{Lemma}
\begin{proof} Let $X$ be a Killing vector field with $(\nabla X)(p)=0$. Then $d\sigma(X)$ and $-X$ are Killing fields having the same value and by (\ref{KillingCalc}) the same derivative at $p$; Proposition \ref{KillingDeterminedByPoint} thus yields $X\in {\mathfrak p}$.

If conversely $X\in {\mathfrak p}$ is given, (\ref{KillingCalc}) immediately shows $(\nabla_Y X)(p)=0$ for all $Y$.  
\end{proof}
\begin{Cor}\label{RSSOnePoint} A Riemannian supermanifold $M$ is symmetric if and only if it is homogeneous and there exists a point $p\in M$ with an isometry $s_p:M\to M$ leaving $p$ fixed and satisfying $d_ps_p=-\id_{T_pM}$.
\end{Cor}

%\begin{Prop} Let $M$ be a Riemannian symmetric superspace. Then the curvature tensor is parallel.
%\end{Prop}
%\begin{proof} Let $p\in M$ and $s_p$ the symmetry at $p$. If $X,Y,Z,W$ are vector fields on $M$, then \begin{align*}(\nabla R)(X,Y,Z,W)(p)&=-d_p s_p (\nabla R)(X,Y,Z,W)(p)\\
%&=-(\nabla R)(ds_p X,ds_p Y,ds_p Z,ds_p W)(p)\\
%&=-(\nabla R)(X,Y,Z,W)(p);\end{align*} thus, the value of $(\nabla R)(X,Y,Z,W)$ vanishes at any point.\marginpar{weiter?}
%\end{proof}

Let $G$ be a connected Lie supergroup and $K$ a closed Lie subsupergroup. Then the pair $(G,K)$ is a {\it symmetric pair} if there exists an involutive automorphism $\sigma$ of $G$ such that $G^{\sigma}_0\subset K\subset G^\sigma$, where $G^\sigma$ is the group of fixed points of $\sigma$.
\begin{Prop} \label{PropRSSViaHomog}Let $(G,K)$ be a symmetric pair and $\left<\cdot,\cdot\right>$ a $G$-invariant graded Riemannian metric on $G/K$. Then $G/K$ is a Riemannian symmetric superspace.
\end{Prop}
\begin{Rem} We do not state that such an invariant metric always exists!
\end{Rem}
\begin{proof} Let $s_K:G/K\to G/K$ be the morphism induced by $\sigma$. More precisely, for $f\in \shO_{G/K}(U)$ we set $s_K^*(f):=\sigma^*(f)$; this is well-defined because of \mbox{$\sigma\circ \Phi=\sigma\circ {\rm pr}_1$}, where the morphisms $\Phi$ and ${\rm pr}_1$ are the same mappings as in (\ref{SheafOfHomSpace}). Then $$s_K\circ \pi=\pi\circ \sigma,$$ where $\pi:G\to G/K$ is the canonical projection. Differentiating this equation at $e\in G$, we get $d_Ks_K=-\id_{T_KG/K}$. In view of Corollary \ref{RSSOnePoint}, it remains to show that $s_K$ is an isometry.

For any $g$ we will write $g:G/K\to G/K$ for left translation with $g$. Let $p=gK_{\red}\in G_{\red}/K_{\red}$, pick vector fields $X$ and $Y$ around $p$ and define \mbox{$X_0:=d{g^{-1}}(X)$} and $Y_0:=d{g^{-1}}(Y)$. If $\rho:G\times G/K\to G/K$ is the standard action, we have \begin{equation}s_K\circ \rho=\rho\circ (\sigma\times s_K)\label{SymSigmaChange}\end{equation} and thus in particular $s_K\circ g={\sigma(g)}\circ s_K$. This yields
\begin{align}\left<ds_K(X_p),ds_K(Y_p)\right>_{s_K(p)}&=\left<ds_K(dg(X_{0,K})),ds_K(dg(Y_{0,K}))\right>_{s_K(p)}\nonumber\\
&=\left<d({\sigma(g)})(ds_K(X_{0,K})),d({\sigma(g)})(ds_K(Y_{0,K}))\right>_{s_K(p)}\nonumber \\
&=\left<d({\sigma(g)})(X_{0,K}),d({\sigma(g)})(Y_{0,K})\right>_{s_K(p)}\nonumber\\
&=\left<X_{0,K},Y_{0,K}\right>_K=\left<X_p,Y_p\right>_{p},\label{SymIsomAtPoint}\end{align} where we used the $G$-invariance of the metric for the last two equalities. 

For any $Z\in {\mathfrak g}$, the vector field $\bar{Z}=(Z_e\otimes I)\circ \rho^*$ is a Killing field because of the $G$-invariance of the metric. Then 
\begin{align*}ds_K(\bar{Z})&=s_K^*\circ (Z_e\otimes I)\circ (s_K\circ \rho)^*=s_K^*\circ (X_e\otimes I)\circ (\sigma\times s_K)^*\circ \rho^*\\
&=(d\sigma(Z)_e\otimes I)\circ \rho^*=\overline{d\sigma(Z)},
\end{align*} so in particular $ds_K(\bar{Z})$ is a Killing field. Thus we may calculate
\begin{align*}&ds_K(\bar{Z})\left<ds_K(X),ds_K(Y)\right>(s_K(p))\\
&\qquad=\left<[ds_K(\bar{Z}),ds_K(X)]_{s_K(p)},ds_K(Y_p)\right>_{s_K(p)}\\
&\qquad\qquad\qquad\qquad+(-1)^{|\bar{Z}||X|}\left<ds_K(X_p),[ds_K(\bar{Z}),ds_K(Y)]_{s_K(p)}\right>_{s_K(p)}\\
&\qquad=\left<[\bar{Z},X]_{p},Y_p\right>_p+(-1)^{|\bar{Z}||X|}\left<X_p,[\bar{Z},Y]_p\right>_p=\bar{Z}\left<X,Y\right>(p),
\end{align*} where we used (\ref{SymIsomAtPoint}). The same calculation shows $$ds_K(\bar{Z_1})\circ\ldots\circ ds_K(\bar{Z_n})\left<ds_K(X),ds_K(Y)\right>(s_K(p))=\bar{Z_1}\circ\ldots\circ\bar{Z_n}\left<X,Y\right>(p)$$ for all $Z_1,\ldots,Z_n\in {\mathfrak g}$. The $G$-invariance of the metric, together with the $G$-homogeneity of $G/K$ shows that there are enough Killing vector fields induced by the action to conclude $s_K^*\left<ds_K(X),ds_K(Y)\right>=\left<X,Y\right>$.
\end{proof}

\subsection{Examples}\label{Section_Examples}

The trivial examples are the following: Clearly, any Riemannian symmetric space is a Riemannian symmetric superspace. Furthermore, $\R^{p|2q}$ with the standard metric is a Riemannian symmetric superspace. Thus, the exterior bundle of the trivial bundle $M\times \R^{2q}\to M$ over a Riemannian symmetric space $M$ gives rise to a Riemannian symmetric superspace since in our language it is merely $M\times \R^{0|2q}$, the product of two such spaces. 

Just like in the standard theory, groups that admit bi-invariant metrics give a class of examples.
\begin{Prop} A Lie supergroup with a bi-invariant graded Riemannian metric is a Riemannian symmetric superspace.
\end{Prop}
\begin{proof} Let $G$ be a Lie supergroup with a bi-invariant graded Riemannian metric. Since the metric is left-invariant, $G$ is homogeneous. The bi-invariance shows that the inverse map $i$ is an isometry (it interchanges left-invariant vector fields with right-invariant vector fields and its differential $di:T_eG\to T_eG$ is minus the identity) and thus gives the symmetry at $e$.
\end{proof}
%In particular, $\SL(n|m)$ for $n,m\ge 1$, $n\neq m$, $\SOSp(n|2m)$ for $n,m\ge 1$, $n\neq 2m+2$,... are Riemannian symmetrc superspaces.
Inspired by the Zirnbauer list \cite{Zirnbauer}, we now present examples of series of Riemannian symmetric superspaces. Note that his notion of Riemannian symmetric superspace is in so far different from ours since he defines them as quotients of {\it complex} Lie supergroups with a certain additional condition. We emphasize that the following list is no attempt of a complete classification -- more examples can for example be given by duality. The classification problem is related to the problem of classifying pseudo-Riemannian symmetric spaces \cite{CahPar}, which is solved in the semi-simple case \cite{Berger} but e.g.~not in the solvable case, see \cite{KathOlbrich}.

\subsubsection{The RSSS $\SL(n|2m)/\SOSp(n|2m)$ ($n\neq 2m$)}

Consider on the Lie superalgebra ${\mathfrak {sl}}(n|2m)$ the involution $\sigma$, given by
$$\left(\begin{array}{c|cc}A & B_1 & B_2 \\ \hline  C_1 & D_1 & D_2 \\ C_2 & D_3 & D_4\end{array}\right)\longmapsto \left(\begin{array}{c|cc}-A^t & C_2^t & -C_1^t \\ \hline -B_2^t & -D_4^t & D_2^t \\ B_1^t & D_3^t &-D_1^t\end{array}\right),$$ where $A, B_i, C_i$ and $D_i$ are $n\times n$-, $n\times m$-, $m\times n$- and $m\times m$-matrices, respectively. An easy calculation shows that $\sigma$ is an automorphism. The decomposition into the eigenspaces of $\sigma$ is $${\mathfrak {sl}}(n|2m)={\mathfrak {osp}}(n|2m)\oplus {\mathfrak p}$$ with ${\mathfrak p}=\{\left(\begin{array}{c|cc}A & B_1 & B_2 \\ \hline B_2^t & D_1 & D_2 \\ -B_1^t & D_3 & D_1^t\end{array}\right)\in {\mathfrak{sl}}(n|2m)\mid A^t=A,\, D_2^t=-D_2,\, D_3^t=-D_3\}.$

The involution $\sigma$ is induced by an isometry of $\SL(n|2m)$: On the level of the underlying Lie groups, the involutive automorphism of $\GL(n)\times \GL(2m)$ given by $$(X,Y)\mapsto ((X^{-1})^t,-J_m(Y^{-1})^tJ_m),$$ %$$\left(\begin{array}{c|c}A & 0 \\ \hline 0 & D\end{array}\right)\mapsto \left(\begin{array}{c|c}(A^{-1})^t & 0 \\ \hline 0 & -J(D^{-1})^t J\end{array}\right),$$
 where $J_m=\left(\begin{array}{cc}0 & I_m \\ -I_m & 0 \end{array}\right)$, restricts to an involutive automorphism of the reduced group $\SL(n|2m)_{\red}$ with $\sigma_0$ as differential and $\SO(n)\times \Sp(m,\R)$ as (connected) fixed point group. Since the question of extending a morphism of Lie superalgebras to a morphism of Lie supergroups concerns only the underlying Lie group (\cite{DelMor}, p. 69), $\sigma$ is induced by an involutive automorphism of  $\SL(n|2m)$ (which is also denoted by $\sigma$). Thus, $(\SL(n|2m),\SOSp(n|2m))$ is a symmetric pair.

In view of Proposition \ref{PropRSSViaHomog} and Theorem \ref{InvMetricHomSpace}, we have to find an $\Ad_{\SL(n|2m)}$-invariant scalar superproduct on $\mathfrak p$. Clearly, the supertrace induces an invariant supersymmetric bilinear form \mbox{$(X,Y)\mapsto \str(XY)$}; we have to show its non-degeneracy. To show the non-degeneracy of its even part, note that the even part of ${\mathfrak{sl}}(n|2m)$ splits as ${\mathfrak{sl}}(n)\oplus{\mathfrak{sl}}(2m)\oplus {\mathfrak u}(1)$. Then $${\mathfrak p}_0={\mathfrak p}_0\cap {\mathfrak{sl}}(n)\oplus {\mathfrak p}_0\cap{\mathfrak{sl}}(2m)\oplus {\mathfrak u}(1)$$ is an orthogonal decomposition with respect to the supertrace. On the first two summands, the supertrace clearly is non-degenerate; on the third, it is so because $n\neq 2m$. To prove the non-degeneracy of the odd part, we define for $$0\neq X=\left(\begin{array}{c|cc}0 & B_1 & B_2 \\ \hline B_2^t & 0 & 0 \\ -B_1^t & 0 & 0\end{array}\right)\in {\mathfrak p}_1$$ an element $$Y=\left(\begin{array}{c|cc}0 & -B_2 & B_1 \\ \hline B_1^t & 0 & 0 \\ B_2^t & 0 & 0\end{array}\right)\in {\mathfrak p}_1$$ to get $\str(XY)=2(\tr(B_1B_1^t)+\tr(B_2B_2^t))>0$. We have thus shown that $\SL(n|2m)/\SOSp(n|2m)$, with the metric induced by the supertrace, is a Riemannian symmetric superspace for $n\neq 2m$. Note that the Killing form is a non-zero multiple of the supertrace since $n\neq 2m$ (see \ref{SectionKilling}) so in this example the Killing form also gives an invariant metric.

\subsubsection{The RSSS $\PSL(2m|2m)/\SOSp(2m|2m)$}
The automorphism $\sigma$, defined as in the previous example, induces an automorphism of $\psl(2m|2m)$. Then the same argumentation as above introduces the structure of a Riemannian symmetric superspace on $\PSL(2m|2m)/\SOSp(2m|2m)$; by passing to the quotient $\PSL(2m|2m)$ of $\SL(2m|2m)$ we achieve that the supertrace is non-degenerate. Note that here we do not have a metric induced by the Killing form since the Killing form of ${\mathfrak{sl}}(2m|2m)$ vanishes, cf. \ref{SectionKilling}.

\subsubsection{The RSSS $\SL(n_1+n_2|m_1+m_2)/\rmS(\GL(n_1|m_1)\times \GL(n_2|m_2))$}
The involution $\sigma$ on the Lie superalgebra ${\mathfrak{sl}}(n_1+n_2|m_1+m_2)$ given by $$\left(\begin{array}{cc|cc}A_1 & A_2 & B_1 & B_2 \\ A_3 & A_4 & B_3 & B_4 \\ \hline  C_1 & C_2 & D_1 & D_2 \\ C_3 & C_4 & D_3 & D_4\end{array}\right)\longmapsto \left(\begin{array}{cc|cc}A_1 & -A_2 & B_1 & -B_2 \\ -A_3 & A_4 & -B_3 & B_4 \\ \hline  C_1 & -C_2 & D_1 & -D_2 \\ -C_3 & C_4 & -D_3 & D_4\end{array}\right)$$ is an automorphism. The corresponding decomposition of ${\mathfrak{sl}}(n_1+n_2|m_1+m_2)$ is
$${\mathfrak{sl}}(n_1+n_2|m_1+m_2)={\mathfrak k}\oplus {\mathfrak p},$$ where \begin{align*}{\mathfrak k}&={\mathfrak s}({\mathfrak{gl}}(n_1|m_1)\times {\mathfrak{gl}}(n_2|m_2))\\ &=\{(X,Y)\in {\mathfrak{gl}}(n_1|m_1)\times {\mathfrak{gl}}(n_2|m_2)\mid \str(X)+\str(Y)=0\}\end{align*} is embedded into ${\mathfrak{sl}}(n_1+n_2|m_1+m_2)$ via $$\left(\left(\begin{array}{c|c}A & B \\ \hline C & D\end{array}\right),\left(\begin{array}{c|c}A' & B' \\ \hline C' & D'\end{array}\right)\right)\mapsto \left(\begin{array}{cc|cc}A & 0 & B & 0 \\ 0 & A' & 0 & B' \\ \hline  C & 0 & D & 0 \\ 0 & C' & 0 & D'\end{array}\right)$$ and the space ${\mathfrak p}$ is given by $${\mathfrak p}=\{\left(\begin{array}{cc|cc}0 & A_1 & 0 & B_1 \\ A_2 & 0 & B_2 & 0 \\ \hline  0 & C_1 & 0 & D_1 \\ C_2 & 0 & D_2 & 0\end{array}\right)\in {\mathfrak{sl}}(n_1+n_2|m_1+m_2)\}.$$ On the level of Lie groups, $$(X,Y)\mapsto (I_{n_1,n_2}XI_{n_1,n_2},I_{m_1,m_2}YI_{m_1,m_2}),$$ where $I_{n,m}=\left(\begin{array}{cc}-I_n & 0 \\ 0 & I_m\end{array}\right)$, is an involution of $\GL(n_1+n_2)\times \GL(m_1+m_2)$ which restricts to an involution of $\SL(n_1+n_2|m_1+m_2)_{\red}$ with fixed point group $\rmS(\GL(n_1)\times \GL(n_2)\times \GL(m_1)\times \GL(m_2))$ and differential $\sigma_0$. Thus, $\sigma$ is induced by an involution of $\SL(n_1+n_2|m_1+m_2)$. Defining the Lie supergroup $\rmS(\GL(n_1|m_1)\times \GL(n_2|m_2))$ to be given by the Harish-Chandra pair $$(\rmS(\GL(n_1)\times \GL(m_1)\times \GL(n_2)\times\GL(m_2)),{\mathfrak s}({\mathfrak{gl}}(n_1|m_1)\times {\mathfrak{gl}}(n_2|m_2))),$$ we thus see that $(\SL(n_1+n_2|m_1+m_2),\rmS(\GL(n_1|m_1)\times \GL(n_2|m_2)))$ is a symmetric pair.

The supertrace $(X,Y)\mapsto \str(XY)$ is an invariant scalar superproduct on ${\mathfrak p}$; the induced invariant metric gives the structure of a Riemannian symmetric superspace. 

Note that for $n_1+n_2=m_1+m_2=:n$, the Killing form vanishes and thus would not be appropriate. In this case, note also that we may write this Riemannian symmetric superspace as $$\PSL(n|n)/{\rm{PS}}(\GL(n_1|m_1)\times \GL(n_2|m_2)),$$ where the ${\rm{P}}$ in ${\rm{PS}}(\GL(n_1|m_1)\times \GL(n_2|m_2))$ means passing to the Lie supergroup obtained by factoring out the one-dimensional center generated by the identity matrix.

\subsubsection{The RSSS $\SOSp(2n|2m)/\rmU(n|m)$}
Write the elements of $\osp(2n|2m)$ in the form $$\left(\begin{array}{cc|cc}A_1 & A_2 & B_1 & B_2 \\ -A_2^t & A_3 & B_3 & B_4 \\ \hline -B_2^t & -B_4^t & C_1 & C_2 \\ B_1^t & B_3^t & C_3 &-C_1^t\end{array}\right),$$ where the $A_i, B_i$ and $C_i$ are $n\times n$-, $n\times m$- and $m\times m$-matrices, respectively, satisfying the relations $A_1^t=-A_1, A_3^t=-A_3, C_2^t=C_2$ and $C_3^t=C_3$. 

The involutive automorphism $\sigma$ of $\osp(2n|2m)$, given by $$\left(\begin{array}{cc|cc}A_1 & A_2 & B_1 & B_2 \\ -A_2^t & A_3 & B_3 & B_4 \\ \hline -B_2^t & -B_4^t & C_1 & C_2 \\ B_1^t & B_3^t & C_3 &-C_1^t\end{array}\right)\mapsto \left(\begin{array}{cc|cc}A_3 & A_2^t & B_4 & -B_3 \\ -A_2 & A_1 & -B_2 & B_1 \\ \hline B_3^t & -B_1^t & -C_1^t & -C_3 \\ B_4^t & -B_2^t & -C_2 & C_1\end{array}\right),$$ yields the decomposition $$\osp(2n|2m)={\mathfrak k}\oplus {\mathfrak p},$$ where $${\mathfrak k}=\{\left(\begin{array}{cc|cc}A_1 & A_2 & B_1 & B_2 \\ -A_2 & A_1 & -B_2 & B_1 \\ \hline -B_2^t & -B_1^t & C_1 & C_2 \\ B_1^t & -B_2^t & -C_2 & C_1\end{array}\right)\mid A_1^t=-A_1,\, A_2^t=A_2,\, C_1^t=-C_1,\, C_2^t=C_2\}$$ and
$${\mathfrak p}=\{\left(\begin{array}{cc|cc}A_1 & A_2 & B_1 & B_2 \\ A_2 & -A_1 & B_2 & -B_1 \\ \hline -B_2^t & B_1^t & C_1 & C_2 \\ B_1^t & B_2^t & C_2 & -C_1\end{array}\right)\mid A_1^t=-A_1,\, A_2^t=-A_2,\, C_1^t=C_1,\, C_2^t=C_2\}$$
Recall that the unitary superalgebra ${\mathfrak u}(n|m)$ is defined as $${\mathfrak u}(n|m)=\{\left(\begin{array}{c|c}A & B \\ \hline -iB^* & C\end{array}\right)\mid A,B,C \text{ complex},\, A^*=-A,\, C^*=-C\}.$$ We can identify it with ${\mathfrak k}$ via $$\left(\begin{array}{c|c}A_1+iA_2 & B_1+iB_2 \\ \hline -B_2^t-iB_1^t & C_1+iC_2\end{array}\right)\mapsto \left(\begin{array}{cc|cc}A_1 & A_2 & B_1 & B_2 \\ -A_2 & A_1 & -B_2 & B_1 \\ \hline -B_2^t & -B_1^t & C_1 & C_2 \\ B_1^t & -B_2^t & -C_2 & C_1 \end{array}\right), $$ where the $A_i,B_i$ and $C_i$ are real matrices. 

On the level of Lie groups, 
$$(X,Y)\mapsto (-J_nXJ_n,-J_mYJ_m)$$ is an automorphism of $\SO(2n)\times \Sp(m;\R)$ with differential $\sigma_0$ and fixed point group isomorphic to $\rmU(n)\times \rmU(m)$. Thus, there is an automorphism of the Lie supergroup $\SOSp(2n|2m)$ turning $(\SOSp(2n|2m),\rmU(n|m))$ into a symmetric pair.

This example somehow plays an extraordinary role since on ${\mathfrak p}_0$, the supertrace is negative definite: $$\str\left(\begin{array}{cc|cc}A_1 & A_2 & 0 & 0 \\ A_2 & -A_1 & 0 & 0 \\ \hline 0 & 0 & C_1 & C_2 \\ 0 & 0 & C_2 & -C_1\end{array}\right)^2=2\, \tr(A_1^2+A_2^2)-2\, \tr(C_1^2+C_2^2).$$ It is non-degenerate on the odd part since for any non-vanishing $$X=\left(\begin{array}{cc|cc}0 & 0 & B_1 & B_2 \\ 0 & 0 & B_2 & -B_1 \\ \hline -B_2^t & B_1^t & 0 & 0 \\ B_1^t & B_2^t & 0 & 0\end{array}\right)\in {\mathfrak p}_1$$ we define $$Y=\left(\begin{array}{cc|cc}0 & 0 & B_2 & -B_1 \\ 0 & 0 & -B_1 & -B_2 \\ \hline B_1^t & B_2^t & 0 & 0 \\ B_2^t & -B_1^t & 0 & 0\end{array}\right)\in {\mathfrak p}_1$$ to get
$$\str(XY)=4\, \tr(B_1B_1^t+B_2B_2^t)>0.$$ Thus, equipped with the metric induced by the negative of the supertrace, $\SOSp(2n|2m)/\rmU(n|m)$ becomes a Riemannian symmetric superspace such that the reduced manifold is a Riemannian symmetric space (the product of one of compact type and one of non-compact type, $\SO(2n)/\rmU(n)\times \Sp(m;\R)/\rmU(m)$).

\subsubsection{The RSSS $\SOSp(n_1+n_2|2m_1+2m_2)/\rmS(\OSp(n_1|2m_1)\times \OSp(n_2|2m_2))$} Writing the elements of $\osp(n_1+n_2|2m_1+2m_2)$ in the 
form 

$$\left(\begin{array}{cc|cccc}A_{11} & A_{12} & B_{11} & B_{12} & B_{13} & B_{14} \\ -A_{12}^t & A_{22} & B_{21} & B_{22} & B_{23} & B_{24} \\ \hline -B_{13}^t & -B_{23}^t & C_{11} & C_{12} & C_{13} & C_{14} \\ -B_{14}^t & -B_{24}^t & C_{21} & C_{22} & C_{14}^t & C_{24}\\ B_{11}^t & B_{21}^t & C_{31} & C_{32} & -C_{11}^t & -C_{21}^t \\ \smash{\underbrace{B_{12}^t}_{n_1}} & \smash{\underbrace{B_{22}^t}_{n_2}} & \smash{\underbrace{C_{32}^t}_{m_1}} & \smash{\underbrace{C_{42}}_{m_2}} & \smash{\underbrace{-C_{12}^t}_{m_3}} & \smash{\underbrace{-C_{22}^t}_{m_4}}\end{array}\right)$$ \vskip12pt 
\noindent with $A_{11}=-A_{11}^t, \, A_{22}=-A_{22}^t, \, C_{13}=C_{13}^t,\, C_{31}=C_{31}^t,\, C_{24}=C_{24}^t$ and \mbox{$C_{42}=C_{42}^t$}, we define an involutive automorphism $\sigma$ of $\osp(n_1+n_2|2m_1+2m_2)$ as follows:
%\begin{tiny}
\begin{align*}&\left(\begin{array}{cc|cccc}A_{11} & A_{12} & B_{11} & B_{12} & B_{13} & B_{14} \\ -A_{12}^t   & A_{22} & B_{21} & B_{22} & B_{23} & B_{24} \\ \hline -B_{13}^t & -B_{23}^t & C_{11} & C_{12} & C_{13} & C_{14} \\ -B_{14}^t & -B_{24}^t & C_{21} & C_{22} & C_{14}^t & C_{24}\\ B_{11}^t & B_{21}^t & C_{31} & C_{32} & -C_{11}^t & -C_{21}^t \\ B_{12}^t & B_{22}^t & C_{32}^t & C_{42} & -C_{12}^t & -C_{22}^t\end{array}\right)\\ &\qquad\qquad\mapsto \left(\begin{array}{cc|cccc}A_{11} & -A_{12} & B_{11} & -B_{12} & B_{13} & -B_{14} \\ A_{12}^t & A_{22} & -B_{21} & B_{22} & -B_{23} & B_{24} \\ \hline -B_{13}^t & B_{23}^t & C_{11} & -C_{12} & C_{13} & -C_{14} \\ B_{14}^t & -B_{24}^t & -C_{21} & C_{22} & -C_{14}^t & C_{24}\\ B_{11}^t & -B_{21}^t & C_{31} & -C_{32} & -C_{11}^t & C_{21}^t \\ -B_{12}^t & B_{22}^t & -C_{32}^t & C_{42}  & C_{12}^t & -C_{22}^t\end{array}\right).
\end{align*} %\end{tiny}
The induced decomposition is $$\osp(n_1+n_2|2m_1+2m_2)=(\osp(n_1|2m_1)\oplus \osp(n_2|2m_2))\oplus {\mathfrak p},$$ where $\osp(n_1|2m_1)\oplus \osp(n_2|2m_2)$ is embedded into $\osp(n_1+n_2|2m_1+2m_2)$ via 
\begin{align*}&\left(\left(\begin{array}{c|cc}A & B_1 & B_2 \\ \hline  -B_2^t & C_1 & C_2 \\ B_1^t & C_3 & -C_1^t\end{array}\right),\left(\begin{array}{c|cc}A' & B'_1 & B'_2 \\ \hline  -{B'}_2^t & C'_1 & C'_2 \\ {B'}_1^t & {C'}_3 & -{C'}_1^t\end{array}\right)\right)\\
&\qquad\qquad\mapsto \left(\begin{array}{cc|cccc}A &0 & B_1 & 0 & B_2 & 0 \\ 0 & A' & 0 & B'_1 & 0 & B'_2 \\ \hline -{B'}_2^t & 0 & C_1 & 0 & C_2 & 0 \\ 0 & -{B'}_2^t & 0 & C'_1 & 0 & C'_2\\ B_1^t & 0 & C_3 & 0 & -C_1^t & 0 \\ 0 & {B_1'}^t & 0 & {C_3'} & 0 & -{C_1'}^t\end{array}\right)
\end{align*} and $${\mathfrak p}=\{\left(\begin{array}{cc|cccc}0 & A_{12} & 0 & B_{12} & 0 & B_{14} \\ -A_{12}^t   & 0 & B_{21} & 0 & B_{23} & 0 \\ \hline 0 & -B_{23}^t & 0 & C_{12} & 0 & C_{14} \\ -B_{14}^t & 0 & C_{21} & 0 & C_{14}^t & 0\\ 0 & B_{21}^t & 0 & C_{32} & 0 & -C_{21}^t \\ B_{12}^t & 0 & C_{32}^t & 0 & -C_{12}^t & 0\end{array}\right)\in \osp(n_1+n_2|2m_1+2m_2)\}.$$ On the level of Lie groups, $$(X,Y)\mapsto (I_{n_1,n_2}XI_{n_1,n_2},L_{m_1,m_2}YL_{m_1,m_2}),$$ where $L_{n,m}=\left(\begin{array}{cccc}-I_n & 0 & 0 & 0 \\ 0 & I_m & 0 & 0 \\ 0 & 0 & -I_n & 0 \\ 0 & 0 & 0 & I_m\end{array}\right)$, is an automorphism of  the reduced Lie group \mbox{$\SO(n_1+n_2)\times \Sp(n;\R)$} with fixed point group $$\rmS(\rmO(n_1)\times\rmO(n_2))\times \Sp(m_1;\R)\times \Sp(m_2;\R)$$ and differential $\sigma_0$. Via the corresponding automorphism of the Lie supergroup $\SOSp(n_1+n_2|2m_1+2m_2)$, the pair $(\SOSp(n_1+n_2|2m_1+2m_2),\rmS(\OSp(n_1|2m_1)\times \OSp(n_2|2m_2))$ becomes a symmetric pair  -- here, $\rmS(\OSp(n_1|2m_1)\times \OSp(n_2|2m_2))$ is the connected component of $\OSp(n_1|2m_1)\times \OSp(n_2|2m_2)$. %is defined via the Harish-Chandra pair $$(\rmS(\rmO(n_1)\times\rmO(n_2))\times \Sp(m_1;\R)\times \Sp(m_2;\R),\osp(n_1|2m_1)\times \osp(n_2|2m_2)).$$
 Again, the supertrace induces an invariant metric on the corresponding homogeneous superspace \mbox{$\SOSp(n_1+n_2|2m_1+2m_2)/\rmS(\OSp(n_1|2m_1)\times \OSp(n_2|2m_2))$} and thus turns it into a Riemannian symmetric superspace.
%\subsection{The Zirnbauer List}
%\subsection{ToDo}
%Einheitliche Bezeichnung: $M\leftrightarrow M_{\red}$, Wert, Auswertung (Varadarajan, p. 141), Halm eindeutig bezeichnen...

\subsubsection{The Exceptional RSSS ${\rm D}(2,1;\alpha)/\SO(2)\times \SOSp(2|2)$}

We also give one example of a family of Riemannian symmetric superspaces not occurring in the Zirnbauer list. The corresponding infinitesimal objects are taken from the tables of Serganova.

The exceptional Lie superalgebra ${\mathfrak g}={\mathfrak d}(2,1;\alpha)$, where $\alpha\in \R\setminus\{0,1\}$ is defined as follows: The even and odd part of ${\mathfrak g}$ are $${\mathfrak g}_0={\mathfrak{sl}}(2)\oplus{\mathfrak{sl}}(2)\oplus{\mathfrak{sl}}(2),$$  $${\mathfrak g}_1=\R^2\otimes\R^2\otimes\R^2,$$ with the ${\mathfrak g}_0$-module structure given by the threefold tensor product of the standard representation of ${\mathfrak{sl}}(2)$ on $\R^2$. The dependence on the parameter $\alpha$ is hidden in the remaining part of the Lie bracket, ${\mathfrak g}_1\times {\mathfrak g}_1\to {\mathfrak g}_0$, see \cite{Scheunert}, Example 5 of Chapter 1, \S 1: Let $\psi:\R^2\times \R^2\to \R$ be the non-degenerate skew-symmetric bilinear form given by $\psi(e_1,e_2)=1$, where $\{e_1,e_2\}$ is the standard basis of $\R^2$. Let $P:\R^2\times \R^2\to {\mathfrak{sl}}(2)$ be the ${\mathfrak{sl}}(2)$-invariant bilinear mapping given by $$P(u,v)w=\psi(v,w)u-\psi(w,u)v$$ for $u,v,w\in \R^2$. Then we define \begin{align*}[u_1\otimes u_2\otimes u_3,v_1\otimes v_2\otimes v_3]=(&\sigma_1 \psi(u_2,v_2)\psi(u_3,v_3)P(u_1,v_1),\\ &\sigma_2 \psi(u_1,v_1)\psi(u_3,v_3)P(u_2,v_2),\\&\sigma_3 \psi(u_1,v_1)\psi(u_2,v_2)P(u_3,v_3)),\end{align*} where the $\sigma_i$ are some real numbers not equal to zero depending on $\alpha$ and satisfying $\sigma_1+\sigma_2+\sigma_3=0$.

We define an involutive automorphism of ${\mathfrak g}$ by $$\sigma=(\tau\oplus\tau\oplus\id_{{\mathfrak{sl}}(2)})\oplus (J_1\otimes J_1\otimes\id_{\R^2}),$$ where $\tau:{\mathfrak{sl}}(2)\to {\mathfrak{sl}}(2)$ is defined by $\tau(A)=-A^t$ and $J_1=\left(\begin{matrix}0 & 1 \\ -1 & 0\end{matrix}\right)$. The ${\mathfrak k}$-part of the corresponding decomposition ${\mathfrak d}(2,1;\alpha)={\mathfrak k}\oplus {\mathfrak p}$ is ${\mathfrak{so}}(2)\oplus \osp(2|2)$.

Let ${\rm D}(2,1;\alpha)$ be the Lie supergroup given by the Harish-Chandra pair $$(\SL(2)\times \SL(2)\times \SL(2),{\mathfrak d}(2,1;\alpha)),$$ where the adjoint representation is the standard one. The automorphism $\sigma$ clearly is induced by an automorphism of ${\rm D}(2,1;\alpha)$ turning $$({\rm D}(2,1;\alpha),\SO(2)\times \SOSp(2|2))$$ into a symmetric pair.

The invariant metric is constructed as follows: Let \mbox{$\left<\cdot,\cdot\right>_1:{\mathfrak g}_1\times {\mathfrak g}_1\to \R$} be the non-degenerate $\ad_{{\mathfrak g}_0}$-invariant skew-symmetric bilinear form given by $$\left<u_1\otimes u_2\otimes u_3,v_1\otimes v_2\otimes v_3\right>_1:=\psi(u_1,v_1)\psi(u_2,v_2)\psi(u_3,v_3);$$ by Proposition \ref{ScalarSP_OddPart}, it extends to an $\ad_{\mathfrak g}$-invariant scalar superproduct $\left<\cdot,\cdot \right>$ on ${\mathfrak g}$. Since $\left<\sigma(u),\sigma(v)\right>_1=\left<u,v\right>_1$ for all $u,v\in {\mathfrak g}_1$ and since $\sigma$ is an automorphism, it follows that $\sigma$ is orthogonal with respect to $\left<\cdot,\cdot\right>$ on the whole of ${\mathfrak g}$. Consequently, ${\mathfrak k}$ is orthogonal to ${\mathfrak p}$ as these spaces are eigenspaces of $\sigma$. Restricting $\sigma$ to ${\mathfrak p}$, we get an $\Ad_{\SO(2)\times \SOSp(2|2)}$-invariant scalar superproduct which induces an invariant metric on ${\rm D}(2,1;\alpha)/\SO(2)\times \SOSp(2|2)$ turning it into a Riemannian symmetric superspace with underlying manifold $\SL(2)/\SO(2)\times \SL(2)/\SO(2)$.

\textsc{Mathematisches Institut, Universit\"at zu K\"oln, Weyertal 86-90, 50931 K\"oln, Germany}

\textsc{Email:} {\texttt{ogoertsc@math.uni-koeln.de}}
%%%%%%%%%%% IM ARTIKEL HIER ZUENDE
\end{document}